\documentclass{amsart}
\usepackage{graphicx} 
\usepackage{amsmath,amssymb,amsthm, mathrsfs}
\usepackage{boxedminipage}
\usepackage[margin=1in]{geometry}
\usepackage{mathpazo}
\usepackage{fancyhdr}
\usepackage[utf8]{inputenc}
\usepackage{color}
\usepackage{tikz}
\usepackage{manfnt}
\usetikzlibrary{positioning, calc}
\usepackage{todonotes}
\usepackage{tabularray}

\newtheorem{theorem}{Theorem}[section]
\newtheorem{proposition}[theorem]{Proposition}
\newtheorem{lemma}[theorem]{Lemma}
\newtheorem{definition}[theorem]{Definition}

\newtheorem{varexample}[theorem]{Example}
\newtheorem{varremark}[theorem]{Remark}
\newtheorem{corollary}[theorem]{Corollary}
\newtheorem{question}[theorem]{Question}

\makeatletter
\def\namedlabel#1#2{\begingroup
	#2%
	\def\@currentlabel{#2}%
	\phantomsection\label{#1}\endgroup
}
\makeatother

\newenvironment{example}{\begin{varexample}
		\begin{normalfont}}{\end{normalfont}
\end{varexample}}

\newenvironment{remark}{\begin{varremark}
		\begin{normalfont}}{\end{normalfont}
\end{varremark}}

\title{RELATING DIFFERENT DEFINITIONS OF LINEAR SERIES ON TROPICAL CURVES}
\author{Eric Burkholder}

\begin{document} 
	\bibliographystyle{alpha}
	
	\begin{abstract}
		We investigate relationships among several recent notions of linear series on tropical curves, including tropical linear series, strongly recursive tropical linear series of Farkas, Jensen, and Payne, and combinatorial limit linear series of Amini and Gierczak. We introduce structured and locally weakly recursive tropical linear series, showing that every locally weakly recursive tropical linear series is a combinatorial limit linear series. Consequently, every strongly recursive tropical linear series is a combinatorial limit linear series. We also obtain an equivalent characterization of combinatorial limit linear series. We establish additional extensions of these results and construct counterexamples showing that the converse implications fail for strongly recursive tropical linear series. Finally, we investigate the extent to which permutation arrays can arise as local combinatorial data of tropical linear series.
	\end{abstract}
	
	\maketitle
	
	\section{Introduction}
	\label{Sec:intro}
	
	Over the past two decades, many researchers have investigated divisors and complete linear series on tropical curves. Their work provided definitions of curves, divisors, and complete linear series on finite and metric graphs that worked analogously to those on algebraic curves. It also provided new proofs of results in classical algebraic geometry through degeneration arguments. A prime example of the former is the tropical Riemann-Roch theorem for finite graphs \cite{baker2007riemann} and for metric graphs \cite{gathmann2008riemann}, \cite{mikhalkin2008tropical}. The tropical proof of the Brill-Noether theorem in \cite{cools2012tropical} was the first significant example of the latter.
	
	Since then, researchers have investigated incomplete linear series on a tropical curve, in analogy to linear series on an algebraic curve. Several definitions of linear series on tropical curves appear in recent literature. One definition of a tropical linear series is provided in \cite{chang2025matroidal}. This definition balances two notions of the rank of a tropical submodule, and all definitions of linear series in the literature are tropical linear series. Meanwhile, Amini and Gierczak in \cite{amini2022limit} defined combinatorial limit linear series. Their definition, given in Subsection \ref{Subsec:CLLS}, focuses on local combinatorial data at points of the metric graph. Both definitions yield finitely generated tropical submodules, whose combinatorial properties are explored in \cite{luo2018idempotent}.
	
	Tropicalizations of linear series satisfy both of these definitions in addition to other properties. Specifically, tropicalizations of linear series are matroidal linear series \cite[Theorem~1.7]{chang2025matroidal}, refined linear series \cite[Theorem~1.7]{amini2022limit}, and strongly recursive tropical linear series \cite[Section~6]{farkas2025kodairadimensionsoverlinemathcalm22overlinemathcalm23}. The first two definitions are also related to generalizations of the theory of matroids: valuated matroids in \cite{chang2025matroidal} and matricube rank functions in \cite{amini2024combinatorial}.
	
	Section \ref{Sec:Prelim} reviews the necessary background on the theory of linear series on tropical curves, as well as the definitions of permutation arrays (adapted from \cite{eriksson2000combinatorial}), slope structures (adapted from \cite{amini2022limit}), tropical linear series, recursive conditions on tropical linear series, structured tropical submodules, and combinatorial limit linear series.
	
	A tropical linear series is \textit{weakly recursive} if any two functions are contained in a tropical linear subseries of rank $1$. It is \textit{locally weakly recursive} if for every point in $\Gamma$ there is a neighborhood around that point such that $\Sigma$ restricted to that neighborhood is weakly recursive. In Section \ref{Sec:LWRTLS-are-structured}, we prove the following result.
	
	\begin{theorem}
		\label{Thm:lwrtls-are-structured}
		Locally weakly recursive tropical linear series are combinatorial limit linear series.
	\end{theorem}
	
	As a special case of Theorem \ref{Thm:lwrtls-are-structured} we arrive at the following.
	
	\begin{corollary}
		\label{Cor:SRTLS-are-structured}
		Strongly recursive tropical linear series are combinatorial limit linear series.
	\end{corollary}
	
	A tropical linear series is \textit{structured} if its slope data and local data are both subject to and store the data of a slope structure. In Section \ref{Sec:Structured-TLS-equiv-to-CLLS}, we prove a key step toward Theorem \ref{Thm:lwrtls-are-structured}.
	
	\begin{theorem}
		\label{Thm:Structured-TLS-equiv-CLLS}
		A tropical linear series is a combinatorial limit linear series if and only if it is structured.
	\end{theorem}
	
	In Section \ref{Sec:extension-and-counterexample} we provide several slight generalizations of Theorem \ref{Thm:lwrtls-are-structured} as well as counterexamples to the converse of Corollary \ref{Cor:SRTLS-are-structured}. In Section \ref{Sec:Realizability} we consider when permutation arrays can be realized as the local data of a tropical linear series.
	
	\begin{figure}
		\begin{center}
			\begin{tikzpicture}[
				every node/.style={
					font=\large,
					align=center
				},
				solid/.style={
					line width=1.1pt
				},
				dashed/.style={
					line width=1.1pt,
					dash pattern=on 6pt off 5pt
				},
				broken/.style={
					line width=1.1pt
				}
				]
				
				
				\node (tls) at (0,7)
				{Tropical Linear Series (TLS)};
				
				\node (clls) at (-3.7,4.5)
				{Combinatorial Limit Linear Series\\
					$=$ Structured TLS};
				
				\node (mls) at (5.0,4.5)
				{Matroidal Linear Series};
				
				\node (refined) at (-6.25,1.15)
				{Refined Linear Series};
				
				\node (lwr) at (0,2.5)
				{Locally Weakly Recursive TLS};
				
				\node (wr) at (0,0.85)
				{Weakly Recursive TLS};
				
				\node (sr) at (0,-0.45)
				{Strongly Recursive TLS};
				
				\node (trop) at (0,-1.9)
				{Tropicalizations of Algebraic Linear Series};
				
				
				\draw[dashed]
				(tls.south west) -- (clls.north);
				
				\draw[dashed]
				(tls.south east) -- (mls.north);
				
				\draw[solid]
				(clls.south) -- (refined.north);
				
				\draw[dashed]
				(clls.south east) -- (lwr.north);
				
				\draw[dashed]
				(lwr.south) -- (wr.north);
				
				\draw[solid]
				(wr.south) -- (sr.north);
				
				\draw[solid] (-0.2,0) -- (0.2,0.2);
				\draw[solid] (-0.2,0.2) -- (0.2,0.4);
				
				\draw[solid]
				(sr.south) -- (trop.north);
				
				\draw[solid]
				(refined.south) -- (trop.north west);
				
				\draw[solid]
				(mls.south) -- (trop.north east);
				
				\draw[solid] (4.05,0.65) -- (4.45,0.85);
				\draw[solid] (4.05,0.85) -- (4.45,1.05);
			\end{tikzpicture}
		\end{center}
		\caption{Diagram Describing Known Relationships between Notions of Linear Series on Tropical Curves}
		\label{Fig:linear-series-poset}
	\end{figure}
	
	Due to the many notions of linear series on tropical curves introduced in this paper and elsewhere in the literature, we have provided Figure \ref{Fig:linear-series-poset} to clarify their known relationships. The diagram is designed to resemble the Hasse diagram of a finite poset. The vertical position of each notion reflects its generality: a notion appearing higher in the diagram is more general than a notion appearing lower in the diagram. A line connecting two notions indicates that one notion implies the other, regardless of whether the line is solid, dashed, or intersected by transverse lines. For example, every tropicalization of an algebraic linear series is a matroidal linear series, and every matroidal linear series is a tropical linear series. By Theorem \ref{Thm:Structured-TLS-equiv-CLLS}, combinatorial limit linear series and structured tropical linear series are equivalent, and are therefore represented by a single element of the diagram.
	
	There are important differences between this diagram and a Hasse diagram. Two notions that appear as incomparable elements in the diagram, such as matroidal linear series and weakly recursive tropical linear series, are not necessarily incomparable as definitions. Rather, their incomparability indicates only that their relationship is not currently characterized by results in the literature. Additionally, notions that appear as distinct elements of the diagram may in fact be equivalent; this possibility is clarified by the additional information encoded in the diagram.
	
	The style of line in the diagram helps record information regarding the status of the relationship between two comparable notions of linear series. A connecting line is dashed, as between tropical linear series and combinatorial limit linear series, when there is a question raised in this paper or elsewhere in the literature as to whether the two notions may in fact be equivalent. A connecting line has two transverse marks, as between tropicalizations of algebraic linear series and matroidal linear series, when the two notions are known not to be equivalent. In this case, there is an example in the literature satisfying one notion but not the other.
	
	We have already identified where some of the relationships in Figure \ref{Fig:linear-series-poset} appear in the literature. Other relationships follow directly from the definitions, while some will be established in this paper. We will refer back to this diagram when introducing new notions of linear series, establishing new relationships, and presenting relevant examples.
	
	\section*{Acknowledgements} This research would not have been possible without the invaluable guidance and advice of David Jensen. The author is especially grateful to the anonymous referee for an exceptionally thorough report. The referee's detailed comments and insightful suggestions led to substantial improvements in the exposition and organization of the manuscript. This work was supported by NSF DMS-2054135.
	
	\section{Preliminaries}
	\label{Sec:Prelim}
	
	\subsection{General Preliminaries}
	\label{Subsec:General-Prelim}
	
	Given a finite graph $G = (V,E)$, a \textit{metric graph} $\Gamma$ with model $G$ is an association of a real interval of finite length to each edge $e \in E$.
	
	A \textit{divisor} on a metric graph is a finitely supported formal sum $$D = \sum_{v \in \Gamma} D(v) \cdot v$$ where $D(v) \in \mathbb{Z}$. A divisor is \textit{effective}, denoted $D \geq 0$, if $D(v) \geq 0$ for all $v \in \Gamma$. The \textit{support} of $D$, denoted $\mathrm{supp}(D)$, is the set of points of $\Gamma$ where $D(v) \neq 0$. The \textit{degree} of $D$ is $$\mathrm{deg}(D) = \sum_{v \in \mathrm{supp}(D)} D(v).$$
	
	For $v \in \Gamma$, let $T_v(\Gamma)$ be the set of outgoing tangent vectors of $v$. Formally, an element $\eta \in T_v(\Gamma)$ is an equivalence class of isometric embeddings $\gamma: [0,\epsilon) \to \Gamma$ satisfying $\gamma(0) = v$. We say an interval $I = [v,w] \subseteq \Gamma$ \textit{contains} a tangent vector $\eta \in T_v(\Gamma)$ if $\gamma([0,\epsilon)) \subset I$ for some representative $\gamma$ of $\eta$. The \textit{valence} of a point $v \in \Gamma$ is $\mathrm{val}(v) = |T_v(\Gamma)|$. For example, a point in the interior of an edge of the underlying finite graph has valence $2$.
	
	Define $\mathrm{PL}(\Gamma)$ to be the set of continuous, piecewise linear functions with integer slopes on $\Gamma$. For $f \in \mathrm{PL}(\Gamma)$ let $\mathrm{sl}_{\eta}(f)$ be the outgoing slope of $f$ along $\eta \in T_v(\Gamma)$. For $v \in \Gamma$, let the \textit{order of vanishing} of $f$ at $v$ be $$\mathrm{ord}_v(f) = - \sum_{\eta \in T_v(\Gamma)} \mathrm{sl}_{\eta}(f).$$ The \textit{principal divisors} on $\Gamma$ are those of the form $$\mathrm{div}(f) = \sum_{v \in \Gamma} \mathrm{ord}_v(f) \cdot v$$ for $f \in \mathrm{PL}(\Gamma)$.
	
	On the set $\mathrm{PL}(\Gamma)$, we define the $||\cdot||_{\infty}$-topology, induced by the norm $$||f||_{\infty} = \max_{v \in \Gamma } \{ |f(v)| \}.$$ Given a divisor $D$ on $\Gamma$, the \textit{complete linear series} is $$R(D) = \{ f \in \mathrm{PL}(\Gamma) : D + \mathrm{div}(f) \geq 0 \}.$$ While $R(D)$ is not a compact set under the $||\cdot||_{\infty}$-topology, if we choose a particular point $v_0 \in \Gamma$, Gathmann and Kerber show in \cite{gathmann2008riemann} that the set $\{ f \in R(D) : f(v_0) = 0 \}$ is compact.
	
	A \textit{tropical linear combination} of functions $f_1, \ldots, f_r \in \mathrm{PL}(\Gamma)$ is $$\theta = \min \{ f_1 + a_1, \ldots , f_r + a_r \}$$ where $a_1, \ldots, a_r$ are real numbers. All sets of functions we consider will be \textit{tropical submodules}: that is, sets closed under tropical linear combinations. Given functions $f_1, \ldots, f_r$, we let $\langle f_1, \ldots, f_r \rangle$ denote the smallest tropical submodule containing these functions.
	
	Analogously, we define tropical notions of independence and dependence. Let $f_1, \ldots, f_r \in \mathrm{PL}(\Gamma)$. If there exist $a_1, \ldots, a_r \in \mathbb{R}$ such that $$\theta = \min\{ f_1 + a_1, \ldots , f_r + a_r \} = \min_{j \neq i} \{ f_j + a_j \}$$ for all $1 \leq i \leq r$, then these functions are \textit{tropically dependent}. That is, $f_1, \ldots, f_r$ are tropically dependent if there exists $a_1, \ldots, a_r$ such that at every point $v \in \Gamma$, $\theta$ is equal to at least two of $f_i + a_i$. Otherwise, the functions are \textit{tropically independent}. In Figure \ref{Fig:Ex-Trop-Lin-Series}, there are three functions on an interval. Any two of these functions are tropically independent, while all three are tropically dependent. We can find the maximum number of tropically independent functions in a finitely generated tropical submodule using the following result.
	
	\begin{lemma}\cite[Lemma~3.2]{chang2025matroidal}
		\label{Lem:generators-determine-dependence}
		Suppose every set of $r+2$ functions in $S = \{ f_1, \ldots, f_s \} \subseteq \mathrm{PL}(\Gamma)$ is tropically dependent. Then every set of $r+2$ functions in $\langle S \rangle$ is tropically dependent.
	\end{lemma}
	
	\subsection{Permutation Arrays}
	\label{Subsec:Permutation-Arrays}
	
	The following subsection follows \cite{eriksson2000combinatorial} closely, with some changes in definitions and conventions that make applications to tropical linear series easier. We formulate our definitions for a rank $r+1$ permutation array. This is to provide some consistency with the rest of the paper, as we will use rank $r+1$ permutation arrays to describe the data of rank $r$ tropical linear series.
	
	Let $[r+1] = \{ 1, \ldots, r+1 \}$. A \textit{$d$-dimensional dot array} $P$ is a subset of $[r_1 + 1] \times [r_2 + 1] \times \cdots \times [r_d + 1]$. We only consider when $r_1 = \cdots = r_d = r$ and write $P \subseteq [r+1]^d$. Given a dot array $P$, we say that an element $\mathbf{x} = (x_1, \ldots, x_d) \in [r+1]^d$ is \textit{dotted} if $\mathbf{x} \in P$ and \textit{empty} if $\mathbf{x} \not \in P$. The set $[r+1]^d$ has a partial order given by $\mathbf{x} \preceq \mathbf{y}$ if $x_i \leq y_i$ for all $1 \leq i \leq d$. The poset's meet operation is pointwise minimum, that is $\mathbf{x} \wedge \mathbf{y} = \mathbf{z}$, where $z_i = \min \{ x_i, y_i \}$.
	
	The \textit{principal subarray} of $P$ at $\mathbf{x} \in [r+1]^d$, denoted $P[\mathbf{x}]$, is the set of elements $\mathbf{y} \in P$ such that $\mathbf{y} \succeq \mathbf{x}$. That is, it is an upper interval of the poset.
	
	Eriksson and Linusson define the \textit{rank of $P$ along the $i$-axis} to be $$\mathrm{rk}_{i}P = \#\{ n : x_i = n \text{ for some element $\mathbf{x} \in P$} \}.$$ That is, $\mathrm{rk}_{i}P$ is the number of distinct values of the $i^{th}$-coordinate of the elements of $P$. Then, $P$ is \textit{rankable of rank} $s+1$, denoted $\mathrm{rank}P = s+1$, if $\mathrm{rk}_{i}P = s+1$ for all $1 \leq i \leq d$.
	
\begin{figure}
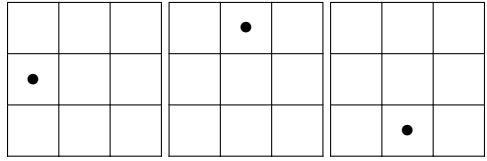

		\centering
		\begin{tblr}{
				rows = {1.5em, rowsep = 2pt},
				columns = {1.5em, colsep = 2pt},
				cells = {m,c},
				hlines,
				vlines,
			}
			&  &  \\
			$\bullet$ &  &  \\
			&  &  \\
		\end{tblr}
		\begin{tblr}{
				rows = {1.5em, rowsep = 2pt},
				columns = {1.5em, colsep = 2pt},
				cells = {m,c},
				hlines,
				vlines,
			}
			& $\bullet$ &  \\
			&  &  \\
			&  &  \\
		\end{tblr}
		\begin{tblr}{
				rows = {1.5em, rowsep = 2pt},
				columns = {1.5em, colsep = 2pt},
				cells = {m,c},
				hlines,
				vlines,
			}
			&  &  \\
			&  &  \\
			& $\bullet$ &  \\
		\end{tblr}
		\caption{Non-Rankable Dot Array $P$}
		\label{Fig:nonrankable-dot-array}
	\end{figure}

	Figure \ref{Fig:nonrankable-dot-array} illustrates an example from \cite[Section~2]{eriksson2000combinatorial} of a non-rankable array. Let the rows be the first dimension, the columns be the second dimension, and the layers be the third dimension. As in the style of \cite{amini2022limit}, we use the convention that the point in the bottom-left corner of the array is the minimum element of the poset $[r+1]^d$, and the point in the top-right corner of the array is the maximum element of the poset $[r+1]^d$. So, the set of dotted points of $P$ in Figure \ref{Fig:nonrankable-dot-array} is $\{ (2,1,1), (3,2,2), (1,2,3) \}$. Then $\mathrm{rk}_{1}P = \mathrm{rk}_3 P = 3$ and $\mathrm{rk}_2 P = 2$, so the dot array is not rankable.
	
	\begin{definition}\cite{eriksson2000combinatorial}
		A dot array is totally rankable if every principal subarray of $P$ is rankable. If $P$ is totally rankable, then the rank array $\rho_P$ of $P$ is the function $\rho_P : [r+1]^d \to \mathbb{Z}$ where $\rho_P(\mathbf{x}) = \mathrm{rank}P[\mathbf{x}]$.
	\end{definition}
	
\begin{figure}
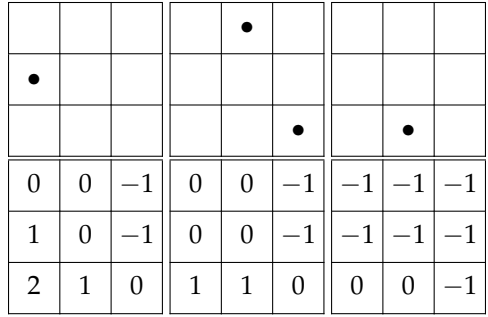

		\centering
		\begin{tblr}{
				rows = {1.5em, rowsep = 2pt},
				columns = {1.5em, colsep = 2pt},
				cells = {m,c},
				hlines,
				vlines,
			}
			&  &  \\
			$\bullet$ &  &  \\
			&  &  \\
		\end{tblr}
		\begin{tblr}{
				rows = {1.5em, rowsep = 2pt},
				columns = {1.5em, colsep = 2pt},
				cells = {m,c},
				hlines,
				vlines,
			}
			& $\bullet$ &  \\
			&  &  \\
			&  & $\bullet$ \\
		\end{tblr}
		\begin{tblr}{
				rows = {1.5em, rowsep = 2pt},
				columns = {1.5em, colsep = 2pt},
				cells = {m,c},
				hlines,
				vlines,
			}
			&  &  \\
			&  &  \\
			& $\bullet$ &  \\
		\end{tblr}\\
		
		\begin{tblr}{
				rows = {1.5em, rowsep = 2pt},
				columns = {1.5em, colsep = 2pt},
				cells = {m,c},
				hlines,
				vlines,
			}
			$1$ & $1$ & $0$ \\
			$2$ & $1$ & $0$ \\
			$3$ & $2$ & $1$ \\
		\end{tblr}
		\begin{tblr}{
				rows = {1.5em, rowsep = 2pt},
				columns = {1.5em, colsep = 2pt},
				cells = {m,c},
				hlines,
				vlines,
			}
			$1$ & $1$ & $0$ \\
			$1$ & $1$ & $0$ \\
			$2$ & $2$ & $1$ \\
		\end{tblr}
		\begin{tblr}{
				rows = {1.5em, rowsep = 2pt},
				columns = {1.5em, colsep = 2pt},
				cells = {m,c},
				hlines,
				vlines,
			}
			$0$ & $0$ & $0$ \\
			$0$ & $0$ & $0$ \\ 
			$1$ & $1$ & $0$ \\ 
		\end{tblr}
		\caption{Totally Rankable Dot Array and Its Rank Array}
		\label{Fig:totally-rankable-dot-array}
	\end{figure}

	In Figure \ref{Fig:totally-rankable-dot-array}, we provide an example of a totally rankable array and its rank array, again taken from \cite[Section~2]{eriksson2000combinatorial}. The totally rankable array in Figure \ref{Fig:totally-rankable-dot-array} is the set of elements $$\{ (2,1,1), (3,2,2), (1,3,2), (1,2,3) \}.$$ The rank array is non-increasing with respect to the poset order, which is also shown in Figure \ref{Fig:totally-rankable-dot-array}.
	
	Eriksson and Linusson provide a characterization of totally rankable dot arrays.
	
	\begin{lemma}\cite[Theorem~3.2]{eriksson2000combinatorial}
		\label{Lem:Totally-Rankable-Characterization}
		A dot array $P$ is totally rankable if and only if, for every two elements $\mathbf{x}, \mathbf{y} \in P$ and two coordinate indices $i$ and $j$ such that $x_i > y_i$ and $x_j = y_j$, there exists $\mathbf{z} \in P$ such that $\mathbf{z} \succeq \mathbf{x} \wedge \mathbf{y}$ and $z_i = y_i$ and $z_j > y_j$.
	\end{lemma}
	
	In order to define permutation arrays, \cite{eriksson2000combinatorial} first define a notion of redundant points of a totally rankable array.
	
	\begin{definition}
		An element $\mathbf{x} \in [r+1]^d$ is redundant if $\mathbf{x} = \bigwedge \mathcal{H}$ for some $\mathcal{H} \subseteq P$ such that $|\mathcal{H}| \geq 2$ and every member in $\mathcal{H}$ has at least one coordinate in common with $\mathbf{x}$. That is, an element is redundant if it can be written as the meet of elements of $P$ in a nontrivial way. For a totally rankable array $P$, we let $R(P)$ denote the redundant positions of $P$.
	\end{definition}
	
	Notice that redundant points need not be elements of $P$. On $[r+1]^d$, Eriksson and Linusson define a \textit{permutation array} of rank $r+1$ and dimension $d$ to be a totally rankable dot array of rank $r+1$ which contains no redundant points.
	
	\begin{lemma}\cite[Proposition~4.1]{eriksson2000combinatorial}
		\label{Lem:Totally-Rankable-Arrays-Have-A-Permutation-Array}
		Two totally rankable dot arrays $P$ and $P'$ have the same rank array if and only if $P \backslash R(P) \subseteq P' \subseteq P \cup R(P)$. In particular, every totally rankable dot array contains a unique permutation array.
	\end{lemma}
	
	By Lemma \ref{Lem:Totally-Rankable-Arrays-Have-A-Permutation-Array}, permutation arrays are the minimal totally rankable dot arrays. For $d=1$, a permutation array is just an array of dimension one with every point dotted. For $d=2$, a permutation array is the set of nonzero elements of a permutation matrix. Figure $\ref{Fig:totally-rankable-dot-array}$ shows an example for $d=3$, as this totally rankable dot array has no redundant points. 
	
\begin{figure}
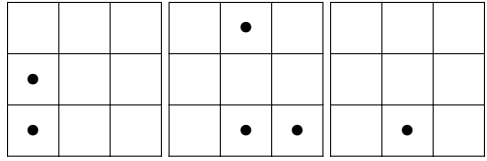

		\centering
		\begin{tblr}{
				rows = {1.5em, rowsep = 2pt},
				columns = {1.5em, colsep = 2pt},
				cells = {m,c},
				hlines,
				vlines,
			}
			&  &  \\
			$\bullet$ &  &  \\
			$\bullet$ &  &  \\
		\end{tblr}
		\begin{tblr}{
				rows = {1.5em, rowsep = 2pt},
				columns = {1.5em, colsep = 2pt},
				cells = {m,c},
				hlines,
				vlines,
			}
			& $\bullet$ &  \\
			&  &  \\
			& $\bullet$ & $\bullet$ \\
		\end{tblr}
		\begin{tblr}{
				rows = {1.5em, rowsep = 2pt},
				columns = {1.5em, colsep = 2pt},
				cells = {m,c},
				hlines,
				vlines,
			}
			&  &  \\
			&  &  \\
			& $\bullet$ &  \\
		\end{tblr}
		\caption{Redundant Closure of Figure \ref{Fig:totally-rankable-dot-array}}
		\label{Fig:redundant-closure}
	\end{figure}

	Eriksson and Linusson also introduce \textit{sparse permutation arrays}, which are arrays with exactly $r+1$ elements. Note that all dimension $2$ permutation arrays are sparse permutation arrays. The \textit{redundant closure} of a permutation array $P$ is $\overline{P} =  P \cup R(P)$. By definition, to obtain the redundant closure of a permutation array $P$, add all of the points of $[r+1]^d$ obtained by applying the meet operation to subsets of $P$. Figure \ref{Fig:redundant-closure} provides the redundant closure of Figure \ref{Fig:totally-rankable-dot-array}, with redundant points at $(1,1,1)$ and $(1,2,2)$.
	
	In \cite{amini2022limit}, Amini and Gierczak define different objects which hold equivalent data to that of a permutation array and its rank array. We rewrite some of the definitions and results of \cite{amini2022limit} in the terminology of \cite{eriksson2000combinatorial}.
	
	\begin{definition}
		The standard rank array of rank $r+1$ and dimension $d$ is the function on $[r+1]^d$ given by $$\rho_{\mathrm{st}}(\mathbf{x}) = \max \{ 0, r + d + 1 - x_1 - x_2 - \cdots - x_{d} \} \text{ for all } \mathbf{x} \in [r+1]^d.$$ The standard permutation array of rank $r+1$ and dimension $d$ is the dot array whose rank array is the standard rank array. That is, the standard permutation array is the dot array $$P = \left \{ \mathbf{x} \in [r+1]^d : \sum_{i=1}^d x_i = r+2 \right \}.$$
	\end{definition}
	
	For $d=2$, the standard permutation arrays are the anti-diagonal matrices, i.e. the matrices with nonzero elements $\{(r+1,1), (r,2), \cdots, (2,r),(1,r+1)\}$.
	
	\begin{lemma}\cite[Proposition~2.14]{amini2022limit}
		\label{Lem:Jumps-Form-Graded-Poset}
		For a permutation array $P$, the set $\overline{P}$ is a graded poset. In particular, if $\mathbf{x} \prec \mathbf{y}$ are two distinct elements of $\overline{P}$, then $\rho_P(\mathbf{x}) > \rho_P(\mathbf{y})$.
	\end{lemma}
	
	\begin{lemma}\cite[Lemma~5.12]{amini2022limit}
		\label{Lem:Jump-of-Same-Rank}
		For a permutation array $P$, let $\mathbf{x}$ be any point in $[r+1]^d$ such that $\rho_P(\mathbf{x}) \geq 1$. Then there exists a unique element $\mathbf{y} \in \overline{P}$ of rank $\rho_P(\mathbf{x})$ with $\mathbf{y} \succeq \mathbf{x}$.
	\end{lemma}
	
	\subsection{Slope Structures}
	\label{Subsec:Slope-Structures}
	
	We use permutation arrays to define slope structures on metric graphs. When Amini and Gierczak define slope structures, they first provide a definition for finite graphs, and then extend that definition to metric graphs. For conciseness, we have written this definition for metric graphs only.
	
	\begin{definition}
		Let $\Gamma$ be a metric graph. Let $V$ be a finite set of points of $\Gamma$ that contains all points of valence not equal to 2. An $(r+1)$-slope structure $\mathfrak{S} = \{ (P_v, S^{\eta_1}, \ldots, S^{\eta_{\mathrm{val}(v)}}) : v \in \Gamma, \eta_{i} \in T_v(\Gamma) \}$ on $(\Gamma,V)$ is the data of
		\begin{enumerate}
			\item [\namedlabel{prop:tangent-vector}{(1)}] For every tangent vector $\eta \in T_v(\Gamma)$, a collection $S^\eta$ of $r+1$ integers $s_j^{\eta}$ such that $$s_1^{\eta} < s_2^{\eta} < \cdots s_{r+1}^{\eta}.$$ 
			\item [\namedlabel{prop:ordering}{(2)}] For every point $v \in \Gamma$, an ordering of the elements of $T_v(\Gamma)$ as $\eta_1, \ldots, \eta_{\mathrm{val}(v)}$.
			\item [\namedlabel{prop:perm-array}{(3)}] For every point $v \in \Gamma$, a rank $r$ and dimension $\mathrm{val}(v)$ permutation array $P_v$.
		\end{enumerate}
		
		This data is subject to the following conditions. 
		
		\begin{itemize}
			\item For every connected component of $\Gamma \backslash V$, and for every tangent vector oriented in the same direction of this component, $S^{\eta}$ is constant. Further, if $\xi$ is a tangent vector in the opposite direction in this component, then $s_j^{\eta} + s_{r+1-j}^{\xi} = 0$.
			\item For every point $v \not \in V$, $P_v$ is the standard rank $r+1$ and dimension $2$ permutation array.
		\end{itemize}
		
	\end{definition}
	
	We include \ref{prop:ordering} only to simplify the definitions that follow. The additional conditions, while lengthy, will not be of particular concern, as they follow naturally when one begins inspecting the local structure of finitely generated tropical submodules.
	
	\begin{example}
		\label{Ex:slope-structure}
		Let $\Gamma$ be the interval metric graph of length $6$ with endpoints $v$ and $u$. We fix this length to provide helpful structure for what will be a running example. Define the following permutation arrays:
		
		$$P^0 = \{ (1), (2) \}$$
		
		$$P^1 = \{ (1,2), (2,1) \}$$
		
		$$P^2 = \{ (1,1), (2,2) \}$$
		
		Let $V = \{ v, u, \overline{3} \}$, where $\overline{3}$ represents the point that is distance $3$ from $v$. Let $\eta_{L}$ denote the left-directed tangent vector based at a point and let $\eta_{R}$ denote the right-directed tangent vector based at a point. We consider the $2$-slope structure with the following description.
		
		At $v$, $P_v = P^0$ and $S^{\eta_R} = \{1,2\}$. For all $w \in (0,3)$, $P_{\overline{w}} = P^1$, $S^{\eta_L} = \{-2,-1\}$, and $S^{\eta_R} = \{1,2\}$. At $\overline{3}$, $P_{\overline{3}} = P^2$, $S^{\eta_L} = \{-2,-1\}$, and $S^{\eta_R} = \{0,1\}$. For all $w \in (3,6)$, $P_{\overline{w}} = P^1$, $S^{\eta_L} = \{-1,0\}$, and $S^{\eta_R} = \{0,1\}$. At $u$, $P_u = P^0$ and $S^{\eta_L} = \{-1,0\}$.
		
		Figure \ref{Fig:Ex-Slope-Structure} provides the most relevant data for defining $\mathfrak{S}$. In this image, the slope lists on each subinterval are assuming a rightward orientation.
		
		On the subinterval $[0,3]$, the vector $S^{\eta_R}$ is the same for every rightward tangent vector and the vector $S^{\eta_L}$ is the same for every leftward tangent vector. Further, $s^R_1 + s^L_2 = 1 + (-1) = 0$ and $s^R_2 + s^L_1 = 2 + (-2) = 0$, which satisfies the first of the slope structure data constraints. The permutation array at the points in $(0,3)$ is $P^1$, which is the standard permutation array. Similar observations are possible on the subinterval $[3,6]$.
		
		The permutation array is standard for all points not in $V$, which satisfies the second of the slope structure data constraints.
	\end{example}
	
	\begin{figure}
		\centering
		\begin{tikzpicture}[scale=2]
			
			\draw[thick] (0,0.5)--(6,0.5);
			
			\coordinate[label=below:$v$] (v) at (0,0.5);
			\coordinate[label=below:$u$] (u) at (6,0.5);
			\coordinate[label=below:$\overline{3}$] (w) at (3,0.5);
			
			\node[circle,fill,inner sep=1pt] at (v) {};
			\node[circle,fill,inner sep=1pt] at (u) {};
			\node[circle,fill,inner sep=1pt] at (w) {};
			
			\node[above=2mm of v] {
				\begin{tblr}{
						rows={1.5em,rowsep=2pt},
						columns={1.5em,colsep=2pt},
						cells={m,c},
						hlines,
						vlines,
					}
					$\bullet$\\
					$\bullet$
				\end{tblr}
			};
			
			\node[above=2mm of u] {
				\begin{tblr}{
						rows={1.5em,rowsep=2pt},
						columns={1.5em,colsep=2pt},
						cells={m,c},
						hlines,
						vlines,
					}
					$\bullet$\\
					$\bullet$
				\end{tblr}
			};
			
			\node[above=2mm of w] {
				\begin{tblr}{
						rows={1.5em,rowsep=2pt},
						columns={1.5em,colsep=2pt},
						cells={m,c},
						hlines,
						vlines,
					}
					& $\bullet$\\
					$\bullet$ &
				\end{tblr}
			};
			
			\node[above=2mm] at (1.5,0.5) {$\{1,2\}$};
			\node[above=2mm] at (4.5,0.5) {$\{0,1\}$};
			
		\end{tikzpicture}
		\caption{Example $2$-Slope Structure on the Interval}
		\label{Fig:Ex-Slope-Structure}
	\end{figure}
	
	A function $f \in \mathrm{PL}(\Gamma)$ is \textit{compatible} with $\mathfrak{S}$ if it satisfies the conditions:
	
	\begin{itemize}
		\item for every point $v \in \Gamma$ and each $\eta \in T_v(\Gamma)$, $\mathrm{sl}_{\eta}(f)$ is an element of $S^{\eta}$.
		
		\noindent Let $\partial_v(f,\eta_i) = j$ when $\mathrm{sl}_{\eta_i}(f) = s_{j}^{\eta_i}$. Define $$\partial_v(f) = (\partial_v(f,\eta_1), \ldots, \partial_v(f,\eta_{\mathrm{val}(v)})).$$
		
		\item for every point $v \in \Gamma$, the vector $\partial_v(f)$ is an element of $\overline{P_v}$.
		
	\end{itemize}
	
	We will use the following language to describe the relationship between a function $f \in \mathrm{PL}(\Gamma)$ and $\partial_v(f) \in \overline{P_v}$. A function $f \in \mathrm{PL}(\Gamma)$ and the vector $\partial_v(f) \in \overline{P_v}$ for some $v \in \Gamma$ are \textit{associated}. Given $\mathbf{x} \in \overline{P_v}$, $f \in \mathrm{PL}(\Gamma)$ \textit{realizes} $\mathbf{x}$ if $\partial_v(f) = \mathbf{x}$. A set $\Sigma \subseteq \mathrm{PL}(\Gamma)$ \textit{realizes} a slope structure $\mathfrak{S}$ if for all $v \in \Gamma$ and for all $\mathbf{x} \in \overline{P_v}$ there exists $f_{v, \mathbf{x}} \in \Sigma$ that realizes $\mathbf{x}$.
	
	If $f \in \mathrm{PL}(\Gamma)$ is compatible with $\mathfrak{S}$, the \textit{$\mathfrak{S}$-rank} of $f$ at $v \in \Gamma$ is the rank of its associated point $\mathbf{x} \in \overline{P_v}$. Therefore, the \textit{$\mathfrak{S}$-rank} of $f$ is at least $1$ for all $v \in \Gamma$. That is, $\rho_v(f) = \rho_v(\partial_v(f))$ where $\rho_v$ is the rank array of $P_v \in \mathfrak{S}$.
	
	\begin{example}
		\label{Ex:compatibility}
		Suppose a function $f$ is compatible with the slope structure defined in Example \ref{Ex:slope-structure}. We impose a few constraints. First, every rightward tangent vector in the interval $[0,3]$, $\mathrm{sl}_{\eta_R}(f)$ is $1$ or $2$ and at every rightward tangent vector in the interval $[3,6]$, $\mathrm{sl}_{\eta_R}(f)$ is $0$ or $1$.
		
		It is easier to describe which functions violate the permutation array condition. At the points in $(0,3)$, the only element of $[2]^2$ not in $\overline{P^1}$ is $(2,2)$. This means that no function can have slope $-1$ on tangent vector $\eta_L$ and slope $2$ on tangent vector $\eta_R$. That is, no function can increase slope moving from left to right. The same condition exists at the points in $(3,6)$.
		
		At the point $\overline{3}$, the elements $(1,1), (2,2) \in P^2$ correspond to functions with slopes $-2$ and then $0$ and functions with slopes $-1$ and then $1$. So, while the functions
		
		\[
		\begin{aligned}
			f_1(x) &=
			\begin{cases}
				x+1, & 0\le x\le 4,\\
				5, & 4\le x\le 6,
			\end{cases}
			\\[1.5em]
			f_2(x) &=
			\begin{cases}
				2x, & 0\le x\le 3,\\
				6, & 3\le x\le 6,
			\end{cases}
		\end{aligned}
		\] 
		are compatible with slope structure $\mathfrak{S}$, the function 
		\[
		\begin{aligned}
			g(x) &=
			\begin{cases}
				x+1, & 0\le x\le 3,\\
				4, & 3\le x\le 6,
			\end{cases}
		\end{aligned}
		\] 
		is not compatible with the slope structure, as $\partial_{\overline{3}}(g) = (2,1) \not \in P^2$.
		
		The $\mathfrak{S}$-rank of $f_1$ is $2$ at $v$ and at the point $\overline{4}$ and is $1$ at all other points in the interval.
	\end{example}
	
	Amini and Gierczak define $$R(\mathfrak{S}) = \{ f \in \mathrm{PL}(\Gamma) : f \text{ compatible with } \mathfrak{S} \},$$ and $R(D, \mathfrak{S}) = R(D) \cap R(\mathfrak{S})$. This will be the set on which we define one notion of a linear series on a tropical curve.
	
	\subsection{Tropical Linear Series}
	\label{Subsec:TLS}
	
	First, we provide the definition of a tropical linear series.
	
	\begin{definition}
		\label{Def:Tropical-Linear-Series}
		A tropical linear series of rank $r$ on a metric graph $\Gamma$ is a pair $(D,\Sigma)$ where $D$ is a divisor and $\Sigma \subseteq R(D)$ is a finitely generated tropical submodule satisfying:
		\begin{enumerate}
			\item for every effective divisor $E$ of degree $r$, there is some $f \in \Sigma$ such that $\mathrm{div}(f) + D \geq E$; and
			\item every set of $r+2$ functions in $\Sigma$ is tropically dependent.
		\end{enumerate}
	\end{definition}
	
	Often we will refer to $\Sigma$ itself as the tropical linear series and assume the context of an appropriate divisor $D$. We provide an example of a rank $1$ tropical linear series originating from \cite[Example~6.10]{farkas2025kodairadimensionsoverlinemathcalm22overlinemathcalm23}, which also continues our running examples from Examples \ref{Ex:slope-structure} and \ref{Ex:compatibility}.
	
	\begin{example}
		\label{Ex:trop-lin-series}
		Let $D = 2v$. Recall $f_1$ and $f_2$ from Example \ref{Ex:compatibility} and define $f_3$ below: 
		
		\[
		\begin{aligned}
			f_3(x) &=
			\begin{cases}
				2x, & 0\le x\le 2,\\
				x+4, & 2\le x\le 6.
			\end{cases}
		\end{aligned}
		\]
		
		Figure \ref{Fig:Ex-Trop-Lin-Series} illustrates these functions. Let $\Sigma = \langle f_1, f_2, f_3 \rangle$.
		
		A tropical combination of $f_1$ and $f_2$ satisfies condition (1) for $E = \overline{w}$ for $w \in [0, 3)$; $f_2$ satisfies condition (1) for $E = \overline{3}$; and a tropical combination of $f_2$ and $f_3$ satisfies condition (1) for $E = \overline{w}$ for $w \in (3, 6]$.
		
		The generators of $\Sigma$ are tropically dependent as illustrated in Figure \ref{Fig:Ex-Trop-Lin-Series} and thus by Lemma \ref{Lem:generators-determine-dependence}, any three functions in $\Sigma$ are tropically dependent.
	\end{example}
	
	\begin{figure}
		\centering
		\begin{tikzpicture}[scale=0.8]
			\draw [thick] (0,0.5)--(6,0.5);
			\draw (0,1.05)--(2,3.05)--(6,5.05);
			\draw [red] (0,1)--(3,4)--(6,4);
			\draw [blue] (0,1.95)--(4,3.95)--(6,3.95);
			\coordinate[label = left:$f_1$] (f1) at (-0.3,1.95);
			\coordinate[label = left:$f_2$] (f2) at (3,4.3);
			\coordinate[label = right:$f_3$] (f3) at (6.3,5.05);
			\coordinate[label = left:$v$] (v) at (0,0.5);
			\coordinate[label = right:$u$] (u) at (6,0.5);
			\node at (v)[circle,fill,inner sep=1pt]{};
			\node at (u)[circle,fill,inner sep=1pt]{};
		\end{tikzpicture}
		\caption{Generators of a Tropical Linear Series of Rank 1 over an Interval}
		\label{Fig:Ex-Trop-Lin-Series}
	\end{figure}
	
	There are two key results we will use that apply to all tropical linear series.
	
	\begin{lemma}\cite[Proposition~8.3]{chang2025matroidal}
		\label{Lem:r+1-slopes}
		Let $\Sigma \subseteq R(D)$ be a tropical linear series of rank $r$. For each tangent vector $\eta$, the set of slopes $\{ \mathrm{sl}_{\eta}(f) : f \in \Sigma \}$ has size exactly $r+1$.
	\end{lemma}
	
	For a tangent vector $\eta$ of $\Gamma$, we denote the set of slopes along $\eta$ as $\mathrm{sl}_{\eta}(\Sigma)$ and order the slopes as follows:
	$$\mathrm{sl}_{\eta}[1] < \mathrm{sl}_{\eta}[2] < \cdots < \mathrm{sl}_{\eta}[r+1].$$
	
	\begin{lemma}\cite[Lemma~5.1]{jensen2022tropical}
		\label{Lem:Finite-slope-vectors}
		Given a tropical linear series $\Sigma \subseteq R(D)$ on a metric graph $\Gamma$, there exists a finite set $V \subset \Gamma$ such that:
		\begin{enumerate}
			\item $V$ contains $\mathrm{supp}(D)$ and all points of valence different from 2, and
			\item $\mathrm{sl}_{\eta}(\Sigma)$ is constant on each oriented edge of $\Gamma \backslash V$.
		\end{enumerate}
	\end{lemma}
	
	\begin{proof}
		Let $\Sigma = \langle f_1, \ldots, f_n \rangle$. Let $V$ consist of the set of points in the support of $D$, the set of points with valence not equal to $2$, and $\mathrm{supp}(\mathrm{div}(f_i))$ for all $i$. On all oriented edges of $\Gamma \backslash V$, the set $\mathrm{sl}_{\eta}(\Sigma)$ is constant, given by the slopes of the generators of $\Sigma$ on these edges.
	\end{proof}
	
	\subsection{Recursive Properties on Tropical Linear Series}
	\label{Subsec:SRTLS}
	
	The following definition of strongly recursive tropical linear series was introduced by Farkas, Jensen, and Payne in \cite{farkas2025kodairadimensionsoverlinemathcalm22overlinemathcalm23}. At the time, they simply called it a tropical linear series, but it has since been renamed in \cite{chang2025matroidal}.
	
	\begin{definition}\cite[Definition~1.19]{chang2025matroidal}
		\label{Def:Strongly-Recursive-Linear-Series}
		A tropical linear series $(D,\Sigma)$ of rank $r$ on a metric graph $\Gamma$ is strongly recursive if it satisfies the additional properties:
		\begin{enumerate}
			\item [(3)] every set of $r$ functions in $\Sigma$ is contained in a strongly recursive tropical linear subseries of rank $r-1$; and
			\item [(4)] if $S_1$ and $S_2$ are subsets of $\Sigma$ of size $r-1$, then there are strongly recursive tropical linear subseries $\Sigma_1$ and $\Sigma_2$ containing $S_1$ and $S_2$ of rank $r-1$ such that $\Sigma_1 \cap \Sigma_2$ contains a strongly recursive tropical linear series of rank $r-2$.
		\end{enumerate}
	\end{definition}
	
	These properties are vacuously true for $r=0$ and trivially true for $r=1$ for all tropical linear series. While \cite{farkas2025kodairadimensionsoverlinemathcalm22overlinemathcalm23} needs (3) and (4) for its applications of tropical linear series, \cite{chang2025matroidal} discusses the issues with these properties, namely that they are exceptionally difficult to prove due to their recursive nature. We will not use (4) in this paper, and only use (3) to provide a counterexample to one of our conjectures in the case of strongly recursive tropical linear series. Instead, the results of this paper will focus on a significantly weakened version of (3) written below. 
	
	First, we need information about metric subgraphs. Given a metric graph $\Gamma$, a \textit{metric subgraph} is a closed and connected subset of $\Gamma$ equipped with the induced metric. For a point $v \in \Gamma$, the star metric subgraph $\mathrm{Star}_{\epsilon}(v)$ is the metric subgraph of $\Gamma$ with edges of length $\epsilon > 0$ whose finite model is a star graph with $\mathrm{val}(v)$ edges.
	
	\begin{definition}
		\label{Def:LWRTLS}
		A tropical linear series $(D,\Sigma)$ is weakly recursive if every set of $2$ functions in $\Sigma$ is contained in a tropical linear subseries of rank $1$. A tropical linear series $(D,\Sigma)$ on a metric graph $\Gamma$ is locally weakly recursive if for all $v \in \Gamma$ there exists a metric subgraph $\Gamma'$ containing the metric subgraph $\mathrm{Star}_{\epsilon}(v)$ for some $\epsilon > 0$ such that $\Sigma|_{\Gamma'}$ is a weakly recursive tropical linear series.
	\end{definition}
	
	Every strongly recursive tropical linear series is weakly recursive, and every weakly recursive tropical linear series is locally weakly recursive. These relationships are reflected in Figure \ref{Fig:linear-series-poset}.
	
	\subsection{Structured Tropical Linear Series and Combinatorial Limit Linear Series}
	\label{Subsec:CLLS}
	
	We now turn our attention to tropical linear series with assumptions related to a slope structure $\mathfrak{S}$.
	
	\begin{definition}
		\label{Def:Structured}
		A tropical submodule $\Sigma \subseteq \mathrm{PL}(\Gamma)$ is \textit{structured} if for some slope structure $\mathfrak{S}$, $\Sigma \subseteq R(\mathfrak{S})$ and $\Sigma$ realizes $\mathfrak{S}$.
	\end{definition}
	
	\begin{example}
		\label{Ex:structured-tls}
		Continuing our running example, the tropical linear series $\Sigma$ defined in Example \ref{Ex:trop-lin-series} is subject to and realizes the slope structure $\mathfrak{S}$ defined in Example \ref{Ex:slope-structure}. Therefore, $\Sigma$ is a structured tropical linear series.
	\end{example}
	
	In \cite{amini2022limit}, the authors assume an implicit slope structure $\mathfrak{S}$ when defining key notions. Recall from Subsection \ref{Subsec:Slope-Structures} that at every point $v \in \Gamma$, the $\mathfrak{S}$-rank of $f$ at $v$ is denoted $\rho_{v}(f)$.
	
	\begin{definition}
		For a divisor $D$ and an effective divisor $E$, we say that $f \in \mathrm{PL}(\Gamma)$ satisfies the Baker-Norine Rank Property for $E$, or \namedlabel{prop:BNRP}{(BNRP)}, if $\mathrm{div}(f) + D \geq E$. We say that $f \in \mathrm{PL}(\Gamma)$ satisfies the Local Rank Property for $E$, or \namedlabel{prop:LRP}{(LRP)}, if $\rho_v(f) \geq E(v) + 1$ for all $v \in \Gamma$.
	\end{definition}
	
	This allows us to define the following structures for tropical submodules of $R(D)$.
	
	\begin{definition}\cite{amini2022limit}
		A tropical submodule $\Sigma \subseteq R(D)$ is admissible of rank $r$ if it is closed under the topology induced by $||\cdot||_{\infty}$, $\Sigma \subseteq R(\mathfrak{S})$ for some $(r+1)$-slope structure $\mathfrak{S}$, and for every effective divisor $E$ on $\Gamma$ of degree $r$, there exists $f \in \Sigma$ satisfying \ref{prop:BNRP} and \ref{prop:LRP} for $E$.
	\end{definition}
	
	We have the following results relevant to properties of admissible submodules.
	
	\begin{lemma}\cite[Proposition~5.5]{amini2022limit}
		\label{Lem:TLS-top-closed}
		All finitely generated tropical submodules are closed under the $||\cdot||_{\infty}$-topology.
	\end{lemma}
	
	\begin{lemma}\cite[Theorem~5.8]{amini2022limit}
		\label{Lem:Admissible-implies-structured}
		Admissible tropical submodules are structured.
	\end{lemma}
	
	Finally, we have the remaining structure defined in \cite{amini2022limit} we will explore in this paper.
	
	\begin{definition}\cite[Definitions~6.2]{amini2022limit}
		\label{Def:CLLS}
		A combinatorial limit linear series of rank $r$ on a metric graph $\Gamma$ is a pair $(D, \Sigma)$ consisting of a divisor $D$ and a finitely generated admissible submodule $\Sigma \subseteq R(D)$ of rank $r$ such that every $r+2$ elements of $\Sigma$ are tropically dependent.
	\end{definition}
	
	\begin{example}
		The tropical linear series $\Sigma$ from Example \ref{Ex:trop-lin-series} is finitely generated, and thus is closed under the $||\cdot||_{\infty}$-topology by Lemma \ref{Lem:TLS-top-closed}. It is subject to the slope structure $\mathfrak{S}$ defined in Example \ref{Ex:slope-structure}. 
		
		The verification in Example \ref{Ex:trop-lin-series} showing that $\Sigma$ satisfies condition (1) of a tropical linear series also provides functions that satisfy both \ref{prop:BNRP} and \ref{prop:LRP}. We verify this for the divisor $E = \overline{3}$. At this point, $\mathrm{sl}_{\eta_L}(f_2) = -2$ and $\mathrm{sl}_{\eta_R}(f_2) = 0$. So, $f_2$ satisfies \ref{prop:BNRP} for $E$. Further, $\rho_{\overline{3}}(\partial_{\overline{3}}(f_2)) = \rho_{\overline{3}}((1,1)) = 2$ where $\rho_{\overline{3}}$ is the rank array for $P^2$, which was defined in Example \ref{Ex:slope-structure}. So, $f_2$ also satisfies \ref{prop:LRP} for $E$. Therefore, $\Sigma$ is an admissible submodule.
		
		We have already verified that any $r+2$ functions of $\Sigma$ are tropically dependent, and so $\Sigma$ is a combinatorial limit linear series. 
	\end{example}
	
	Combinatorial limit linear series are tropical linear series by definition and structured by Lemma \ref{Lem:Admissible-implies-structured}. The main results of this paper show that tropical linear series with light recursive conditions are structured tropical linear series and structured tropical linear series are combinatorial limit linear series.
	
	\section{Structured Tropical Linear Series are Equivalent to Combinatorial Limit Linear Series}
	\label{Sec:Structured-TLS-equiv-to-CLLS}
	
	In this section we show that all structured tropical linear series are combinatorial limit linear series. To do so, we need to investigate the relationship between a function satisfying the Baker-Norine rank property and the local rank property.
	
	\begin{lemma}
		\label{Lem:LRP-implies-BNRP}
		Let $\Sigma$ be a rank $r$ structured tropical linear series subject to $(r+1)$-slope structure $\mathfrak{S}$ and let $E$ be an effective divisor with $\mathrm{deg}(E) \leq r$. If $f \in \Sigma$ satisfies \ref{prop:LRP} for $E$, then $f$ also satisfies \ref{prop:BNRP} for $E$. 
	\end{lemma}
	
	\begin{proof}
		If $f$ satisfies \ref{prop:LRP} for $E$, then for all $v \in \Gamma$, $\rho_v(f) \geq E(v) + 1$. Let $P_v \in \mathfrak{S}$ be the permutation array at $v \in \Gamma$, let $\rho_{P_v}$ be the rank array of $P_v$. We prove that if $\rho_v(f) \geq k+1$, then $[\mathrm{div}(f) + D](v) \geq k$.
		
		The case $k = 0$ follows because $f \in R(D)$. Now, assume for induction that if $\rho_v(f) \geq k+1$, then $[\mathrm{div}(f) + D](v) \geq k$.
		
		Let $\rho_v(f) = k + 1 \geq 2$. The function $f$ is associated to some $\mathbf{x} \in \overline{P_v}$, and $\rho_v(\mathbf{x}) = k + 1 \geq 2$. By Lemma \ref{Lem:Jumps-Form-Graded-Poset} there must be some $\mathbf{y} \in \overline{P_v}$ such that $\mathbf{y} \succ \mathbf{x}$ and $\rho_v(\mathbf{y}) = \rho_v(\mathbf{x}) - 1 = k$. Because $\Sigma$ realizes the slope structure, there is some $g \in \Sigma$ associated to $\mathbf{y}$. We have $[\mathrm{div}(f) + D](v) > [\mathrm{div}(g) + D](v)$ because $\mathrm{sl}_{\eta_i}(f) < \mathrm{sl}_{\eta_i}(g)$ for some index $i$ and $\mathrm{sl}_{\eta_j}(f) \leq \mathrm{sl}_{\eta_j}(g)$ for all other indices $j$. We know that $[\mathrm{div}(g) + D](v) \geq \rho_v(\mathbf{y}) = k$ by the inductive hypothesis. Therefore, $[\mathrm{div}(f) + D](v) \geq \rho_v(\mathbf{y}) + 1 = k + 1$.
	\end{proof}
	
	The converse of this statement is false. In the example below, we provide a function which satisfies the Baker-Norine rank property but fails the local rank property.
	
	\begin{example}
		\label{Ex:bnrp-and-lrp}
		Let $\Gamma$ be the metric graph with two edges of equal length between $u$ and $v$, as in Figure \ref{Fig:bnrp-and-lrp}. Let $\Sigma \subseteq R(3u)$ be the tropical linear series generated by the constant function and the function with slope $1$ on both edges. The rank array at the point $u$ is in Figure \ref{Fig:bnrp-and-lrp}.
		
		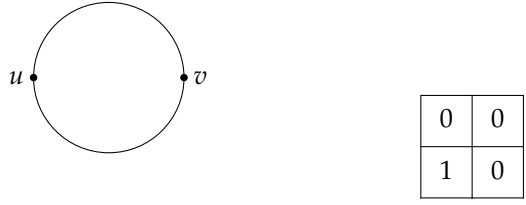
\begin{figure}
			\centering
			
			\begin{minipage}[c]{0.35\textwidth}
				\centering
				\begin{tikzpicture}
					\draw (0,0) circle (1);
					\coordinate[label=left:$u$] (u) at (-1,0);
					\node[circle,fill,inner sep=1pt] at (u) {};
					\coordinate[label=right:$v$] (v) at (1,0);
					\node[circle,fill,inner sep=1pt] at (v) {};
				\end{tikzpicture}
			\end{minipage}
			\begin{minipage}[c]{0.2\textwidth}
				\centering
				\begin{tblr}{
						rows={1.5em, rowsep=2pt},
						columns={1.5em, colsep=2pt},
						cells={m,c},
						hlines,
						vlines,
					}
					$1$ & $1$ \\
					$2$ & $1$ \\
				\end{tblr}
			\end{minipage}
			
			\caption{Graph $\Gamma$ and local rank array at $u$ for Example \ref{Ex:bnrp-and-lrp}}
			\label{Fig:bnrp-and-lrp}
		\end{figure}
		
		The function with slope $1$ on both edges satisfies \ref{prop:BNRP} for $E = u$, but fails \ref{prop:LRP} for $E = u$. Meanwhile, the constant function satisfies both \ref{prop:BNRP} and \ref{prop:LRP} for $E = u$.
	\end{example}
	
	Fortunately, we will be able to prove the existence of a function satisfying the local rank property by taking limits of functions that we know exist in $\Sigma$. We must first prove lemmas about these limits.
	
	\begin{lemma}
		\label{Lem:Limits-1}
		Let $(f_n)_{n \in \mathbb{N}}$ be a sequence of functions in a finitely generated tropical submodule $\Sigma$ on a metric graph $\Gamma$. If there exists $v \in \Gamma$ such that $f_n(v) = 0$ and there exists a slope structure $\mathfrak{S}$ such that $f_n$ is compatible with $\mathfrak{S}$ for all $n$, then there exists a convergent subsequence $(f_m)_{m \in \mathbb{N}}$ such that $\displaystyle \lim_{m \to \infty} f_m = f \in \Sigma$ and $\partial_v(f_m)$ is constant for all $m$.
	\end{lemma}
	
	\begin{proof}
		As mentioned in Subsection \ref{Subsec:General-Prelim}, we know that the set $\{ f \in R(D) : f(v) = 0 \}$ is compact, so the sequence contains a convergent subsequence $(f_m)_{m \in \mathbb{N}}$. By Lemma \ref{Lem:TLS-top-closed}, $\Sigma$ is topologically closed because it is finitely generated, so $\displaystyle \lim_{m \to \infty} f_m = f \in \Sigma$. Because $\overline{P_v}$ is a finite set, we may choose $f_m$ so that $\partial_v(f_m)$ is constant for all $m$.
	\end{proof}
	
	\begin{lemma}
		\label{Lem:Limits-2}
		Let $D$ be an effective divisor on a metric graph $\Gamma$, $(f_n)_{n \in \mathbb{N}} \in R(D)$, $p$ a point in $\Gamma$, and $\eta$ a tangent vector at $p$. If $\displaystyle \lim_{n \to \infty} f_n = f$ and $\displaystyle \lim_{n \to \infty} \mathrm{sl}_{\eta}(f_n)$ exists, then $\displaystyle \mathrm{sl}_{\eta}(f) \leq \lim_{n \to \infty} \mathrm{sl}_{\eta}(f_n)$.
		
		Moreover, if $(p_n)_{n \in \mathbb{N}} \in \Gamma$ where $\displaystyle \lim_{n \to \infty} p_n = p$, all $p_n$ are on $I \backslash \{ p \} = (p,p']$ where $I = [p,p']$ is an interval of the metric graph which contains $\eta$, and $\mathrm{div}(f_n) + D \geq p_n$, then $\displaystyle \mathrm{sl}_{\eta}(f) < \lim_{n \to \infty} \mathrm{sl}_{\eta}(f_n)$.
	\end{lemma}
	
	\begin{proof}
		Let $I = [p,p']$ be the interval of length $\ell > 0$ of the metric graph which contains $\eta$. Let $\epsilon_1$ be the distance along $I$ from $p$ to the nearest point that is not valence-2 or is in the support of $D$. Let $\epsilon_2$ be the distance along $I$ from $p$ to the nearest point where $f$ is not differentiable. Let $\epsilon = \frac{\min \{ \epsilon_1, \epsilon_2, \ell \}}{4}$. Assume for contradiction that $\displaystyle \mathrm{sl}_{\eta}(f) > \lim_{n \to \infty} \mathrm{sl}_{\eta}(f_n)$. 
		
		There exists $N_1 \gg 1$ such that for all $n \geq N_1$, $|f_n(x) - f(x)| < \epsilon$ for all $x \in \Gamma$. Further, since $\displaystyle \mathrm{sl}_{\eta}(f) > \lim_{n \to \infty} \mathrm{sl}_{\eta}(f_n)$, there exists $N \geq N_1$ such that $\mathrm{sl}_{\eta}(f) > \mathrm{sl}_{\eta}(f_{N})$ and $|f_{N}(p) - f(p)| < \epsilon$.
		
		Let $x_0$ be the point along $I$ that is $2\epsilon$ from $p$. Between $p$ and $x_0$ the slope of $f$ cannot decrease by assumption and the slope of $f_N$ cannot increase since no point in this interval is in the support of $D$. Thus, $\mathrm{sl}_{\eta}(f) - \mathrm{sl}_{\eta}(f_{N}) \geq 1$ on the interval $[p,x_0]$. Therefore, 
		
		\begin{align*}
			|f_{N}(x_0) - f(x_0)| & \geq (\mathrm{sl}_{\eta}(f) - \mathrm{sl}_{\eta}(f_{N})) \cdot 2\epsilon - |f_{N}(p) - f(p)| \\
			& \geq 2\epsilon - \epsilon \\
			& = \epsilon.
		\end{align*}
		This is a contradiction. Thus, $\displaystyle \mathrm{sl}_{\eta}(f) \leq \lim_{n \to \infty} \mathrm{sl}_{\eta}(f_n)$.
		
		Assume that $(p_n)_{n \in \mathbb{N}} \in \Gamma$ where $\displaystyle \lim_{n \to \infty} p_n = p$, all $p_n$ are on $I \backslash \{ p \} = (p,p']$, and $\mathrm{div}(f_n) + D \geq p_n$. Assume for contradiction that $\displaystyle \mathrm{sl}_{\eta}(f) = \lim_{n \to \infty} \mathrm{sl}_{\eta}(f_n)$.
		
		There exists $N_2 \gg 1$ such that for all $n \geq N_2$, the distance from $p$ to $p_n$ along $I$ is less than $\epsilon$. Let $N' = \max \{ N_1, N_2 \}$. Let $x_1$ be the point that is $4\epsilon$ away from $p$ along $I$. By the assumption that $\displaystyle \mathrm{sl}_{\eta}(f) = \lim_{n \to \infty} \mathrm{sl}_{\eta}(f_n)$, the fact that no point in the support of $D$ in the interval $[p,p_{N'}]$, and the fact that $f$ is differentiable on $(p,p_{N'}]$, the slope of $f$ is greater than or equal to the slope of $f_{N'}$ on the entire interval $[p,p_{N'}]$. Let $\eta' \in T_{p_{N'}}(\Gamma)$ be the tangent vector at $p_{N'}$ which is contained in the interval $[p_{N'},x_1]$. Because $\mathrm{div}(f_{N'}) + D \geq p_{N'}$ and the fact that $f$ is differentiable on $[p_{N'},x_1)$, $\mathrm{sl}_{\eta'}(f_{N'}) < \mathrm{sl}_{\eta'}(f)$.
		
		On the interval $[p_{N'}, x_1]$, the slope of $f$ cannot decrease by assumption and the slope of $f_{N'}$ cannot increase since no point in this interval is in the support of $D$. Thus, $\mathrm{sl}(f) - \mathrm{sl}(f_{N'}) \geq 1$ on $[p_{N'},x_1]$. Therefore,
		
		\begin{align*}
			|f_{N'}(x_1) - f(x_1)| & \geq (\mathrm{sl}_{\eta}(f) - \mathrm{sl}_{\eta}(f_{N'})) \cdot 3 \epsilon - |f_{N'}(p_{N'}) - f(p_{N'})| \\
			& \geq 3 \epsilon - \epsilon \\
			& = 2 \epsilon
		\end{align*}
		This is a contradiction. Thus, $\displaystyle \mathrm{sl}_{\eta}(f) < \lim_{n \to \infty} \mathrm{sl}_{\eta}(f_n)$.
	\end{proof}
	
	\begin{remark}
		We assumed in Lemma \ref{Lem:Limits-2} that $\displaystyle \lim_{n \to \infty} \mathrm{sl}_{\eta}(f)$ exists. During a discussion of this paper, Omid Amini and Lucas Gierczak observed that if this assumption is dropped, we may modify the lemma's statements to conclude that $\displaystyle \mathrm{sl}_{\eta}(f) \leq \liminf_{n \to \infty} \mathrm{sl}_{\eta}(f_n)$ and $\displaystyle \mathrm{sl}_{\eta}(f) < \liminf_{n \to \infty} \mathrm{sl}_{\eta}(f_n)$, respectively, without requiring any changes to the proofs.
	\end{remark}
	
	\begin{proposition}
		\label{Prop:Structured-implies-BNRP-and-LRP}
		Let $\Sigma$ be a rank $r$ structured tropical linear series. Then for all effective divisors $E$ with $\mathrm{deg}(E) = r$, there exists a function $f \in \Sigma$ satisfying \ref{prop:BNRP} and \ref{prop:LRP} for $E$.
	\end{proposition}
	
	\begin{proof}
		We make the following choices before proving the result. Let $p_1, \ldots, p_r$ be $r$ not necessarily distinct points in $\Gamma$. For each $i = 1, \ldots, r$, choose a tangent vector $\eta_i \in T_{p_i}(\Gamma)$. Let $I_i = [p_i, p_i']$ be an interval which contains $\eta_i$ and let $(p_{i,n})_{n \in \mathbb{N}}$ be a collection of points in $I_i \backslash \{ p_{i} \}$ such that $\displaystyle \lim_{n \to \infty} p_{i,n} = p_i$.
		
		Let $0 \leq k \leq r$. We will prove, by induction on $k$, that there exists $f_{k} \in \Sigma$ satisfying \ref{prop:BNRP} for $E = p_1 + \cdots + p_r$ and \ref{prop:LRP} for $E_{k} = p_1 + \cdots + p_{k}$.
		
		For $k=0$, we know that there exists $f_0 \in \Sigma$ satisfying \ref{prop:BNRP} for $E$ because $\Sigma$ is a tropical linear series. Further, $\rho_v(f_0) \geq 1$ for all $v \in \Gamma$ as this is true for any function compatible with the underlying slope structure, so $f_0$ satisfies \ref{prop:LRP} for $E_0 = 0$.
		
		By the inductive hypothesis, for each $p_{k,n}$ there exists a function $f_{k,n}$ satisfying \ref{prop:BNRP} for $E_{k,n} = E_{k} + p_{k+1,n} + p_{k+2} + \cdots + p_r$ and \ref{prop:LRP} for $E_{k}$. By tropical scaling, we may further require that $f_{k,n}(p_k) = 0$ for all $n$.
		
		By Lemma \ref{Lem:Limits-1}, we may choose a convergent subsequence $(f_{k+1,m})_{m \in \mathbb{N}}$ such that $\displaystyle \lim_{m \to \infty} f_{k+1,m} = f_{k+1} \in \Sigma$ and $\partial_{p_{k+1}}(f_{k+1,m}) \in \overline{P_{p_{k+1}}}$ is constant for all $m$. We now prove several facts about $f_{k+1}$.
		
		First, let $v$ be any point in $\Gamma$. As there are finitely many elements in $\overline{P_v}$, there is an infinite subsequence $(f_{k+1,m_v})_{m_v \in \mathbb{N}}$ such that $\partial_{v}(f_{k+1,m_v}) \in \overline{P_v}$ is constant. Furthermore, this subsequence has the same limit as the original sequence. Therefore for all $\eta \in T_v(\Gamma)$, $\displaystyle \lim_{m_v \to \infty} \mathrm{sl}_{\eta}(f_{k+1,m_v})$ exists, and by Lemma \ref{Lem:Limits-2}, $\displaystyle \mathrm{sl}_{\eta}(f_{k+1}) \leq \lim_{m_v \to \infty} \mathrm{sl}_{\eta}(f_{k+1,m_v})$. In fact, by choosing the subsequence so that $\partial_{v}(f_{k+1,m_v})$ is constant, the sequence $(\mathrm{sl}_{\eta}(f_{k+1,m_v}))_{m \in\mathbb{N}}$ is constant for every tangent vector at $v$. Hence, we may compare $f_{k+1}$ with any element of this subsequence. Let $\mathbf{x}_v, \mathbf{y}_v \in \overline{P_v}$ be associated points of $f_{k+1}, f_{k+1,m_v'}$, respectively. The slope is bounded above on all tangent vectors of $v$, so $\mathbf{x}_v \preceq \mathbf{y}_v$, which implies that $\rho_v(\mathbf{x}_v) \geq \rho_v(\mathbf{y}_v)$ because rank arrays are non-increasing.
		
		Next, consider the case when $v = p_{k+1}$. The argument for any point $v$ of $\Gamma$ given in the previous paragraph applies. Further, by the second part of Lemma \ref{Lem:Limits-2}, $\mathrm{sl}_{\eta_{k+1}} (f_{k+1}) < \mathrm{sl}_{\eta_{k+1}} (f_{k+1,m_{p_{k+1}}'})$, which implies that $\mathbf{x}_{p_{k+1}} \prec \mathbf{y}_{p_{k+1}}$. Therefore, we have the stronger statement that $\rho_{p_{k+1}}(\mathbf{x}_{p_{k+1}}) > \rho_{p_{k+1}}(\mathbf{y}_{p_{k+1}})$ by Lemma \ref{Lem:Jumps-Form-Graded-Poset}.
		
		The above arguments show that $f_{k+1}$ satisfies \ref{prop:LRP} for $E_{k+1}$ and thus also satisfies \ref{prop:BNRP} for $E_{k+1}$ by Lemma \ref{Lem:LRP-implies-BNRP}. It is left to verify that $f_{k+1}$ satisfies \ref{prop:BNRP} for $E$.
		
		We have already shown that at all points $v \in \Gamma$ and all tangent vectors $\eta \in T_v(\Gamma)$, $\mathrm{sl}_{\eta}(f_{k+1}) \leq \mathrm{sl}_{\eta}(f_{k+1,m_v'})$. Therefore, $f_{k+1,m_v'}$ satisfying \ref{prop:BNRP} for $E - p_{k+1}$ implies that $f_{k+1}$ satisfies \ref{prop:BNRP} for $E - p_{k+1}$.
		
		If $p_{k+1}$ is distinct from the points $p_1, \ldots, p_{k}, p_{k+2}, \ldots, p_{r}$, then $f_{k+1}$ satisfies \ref{prop:BNRP} for $E$ as it is already shown to satisfy \ref{prop:BNRP} for $E_{k+1}$ above. Otherwise, $\mathrm{sl}_{\eta_{k+1}}(f_{k+1}) < \mathrm{sl}_{\eta_{k+1}}(f_{k+1,m_{p_{k+1}}'})$ implies that $[\mathrm{div}(f_{k+1})](p_{k+1}) > [\mathrm{div}(f_{k+1,m_{p_{k+1}}'})](p_{k+1})$, which means that $f_{k+1}$ satisfies \ref{prop:BNRP} for $E$.
	\end{proof}
	
	So, we come to the proof of Theorem \ref{Thm:Structured-TLS-equiv-CLLS}.
	
	\begin{proof}[Proof of Theorem~\ref{Thm:Structured-TLS-equiv-CLLS}]
		As previously stated, combinatorial limit linear series are tropical linear series by definition. They are assumed to be admissible submodules, and thus are structured by Lemma \ref{Lem:Admissible-implies-structured}.
		
		A rank $r$ structured tropical linear series is finitely generated by assumption, closed under the topology induced by $||\cdot||_{\infty}$ by Lemma \ref{Lem:TLS-top-closed}, and subject to an $(r+1)$-slope structure by definition. For all effective divisors of degree $r$, there is a function in a structured tropical linear series satisfying both \ref{prop:BNRP} and \ref{prop:LRP} for $E$ by Proposition \ref{Prop:Structured-implies-BNRP-and-LRP}. Any $r+2$ functions of a structured tropical linear series are tropically dependent by definition, and so structured tropical linear series are combinatorial limit linear series.
	\end{proof}
	
	\section{Locally Weakly Recursive Tropical Linear Series are Structured}
	\label{Sec:LWRTLS-are-structured}
	
	In this section we show that all locally weakly recursive tropical linear series are structured, which is then reflected in the diagram in Figure \ref{Fig:linear-series-poset}. To do this, we will recycle some notation and terminology from the definition of an $(r+1)$-slope structure. In particular, we will construct a dot array $\overline{P_v}$ of the local data of certain tropical linear series and prove that it is the redundant closure of a permutation array. We first consider an arbitrary rank $r$ tropical linear series.
	
	Let $\Gamma$ be a metric graph, let $v$ be a valence-$d$ point on $\Gamma$, and let $\eta_1, \ldots, \eta_d$ be the distinct tangent vectors based at $v$. By Lemma \ref{Lem:r+1-slopes}, we may store the local data of a function $f \in \Sigma$ as an element of $[r+1]^d$ in the following way.
	
	Let $\partial_v(f,\eta_i) = j$ when $\mathrm{sl}_{\eta_i}(f) = \mathrm{sl}_{\eta_i}[j]$. Define $\partial_v(f) = (\partial_v(f,\eta_1), \ldots, \partial_v(f,\eta_d))$. Then, we call the \textit{local array} of $\Sigma$ at $v$ the dot array $$\overline{P_v} = \{ \partial_v(f) : f \in \Sigma \}.$$
	
	Again, we say that $f \in \Sigma$ and $\mathbf{x} \in \overline{P_v}$ are \textit{associated} if $\partial_{v}(f) = \mathbf{x}$.
	
	\begin{example}
		Recall the tropical linear series $\Sigma$ defined in Example \ref{Ex:trop-lin-series}. We concluded in Example \ref{Ex:structured-tls} that $\Sigma$ was subject to the slope structure $\mathfrak{S}$ defined in Example \ref{Ex:slope-structure}. Now, we will use the local arrays and slope data of $\Sigma$ to construct a slope structure that $\Sigma$ will necessarily be subject to and realize.
		
		By Lemma \ref{Lem:r+1-slopes}, we know that we may define the tangent vector data using the slopes of the generators of $\Sigma$. The rightward directed slope vector of a point in the interval $[v,3)$ is the set $\{ 1,2 \}$ and the leftward directed slope vector of a point in the interval $(v,3]$ is the set $\{ -2, -1 \}$. The rightward directed slope vector of a point in the interval $[3,u)$ is the set $\{ 0,1 \}$ and the leftward directed slope vector of a point in the interval $(3,u]$ is the set $\{ -1, 0 \}$. 
		
		We determine the local arrays at every point of the interval. As before, we use the notation that the bottom left corner of the array is the minimal element of $[r+1]^d$, the rows are the first dimension and the columns are the second dimension.
		
		At the endpoints $v$ and $u$, we consider a $[2]^1$-array, and there are functions with both possible slopes on the tangent vectors of these points. So, the local array for both of these points is the topmost array in Figure \ref{Fig:rk2-d12-local-arrays}, which is a permutation array.
		
		For all points $\overline{w}$ for $w \in (0,3) \cup (3,6)$, there are functions in $\Sigma$ with slope $sl_L[1]$ on the left tangent vector and slope $sl_R[2]$ on the right tangent vector, functions with slope $sl_L[1]$ on the left tangent vector and slope $sl_R[1]$ on the right tangent vector, and functions with slope $sl_L[2]$ on the left tangent vector and slope $sl_R[1]$ on the right tangent vector. Therefore, the local array of these points is given as the left dot array in Figure \ref{Fig:rk2-d12-local-arrays}. This dot array is the redundant closure of the standard permutation array of rank $2$ and dimension $2$.
		
		At the point $\overline{3}$, the outgoing slopes on the left tangent vector are $-2$ and $-1$, and the outgoing slopes on the right tangent vector are $0$ and $1$. All functions of $\Sigma$ will either locally have slope $-2$ and then slope $0$, or slope $-1$ and then slope $1$. So, the local array at $\overline{3}$ must be the right dot array in Figure \ref{Fig:rk2-d12-local-arrays}, which is a permutation array.
	\end{example}
	
	\begin{figure}
		\centering
		\begin{tblr}{
				rows = {1.5em, rowsep = 2pt},
				columns = {1.5em, colsep = 2pt},
				cells = {m,c},
				hlines,
				vlines,
			}
			$\bullet$ \\
			$\bullet$
		\end{tblr} \\
		
		\begin{tblr}{
				rows = {1.5em, rowsep = 2pt},
				columns = {1.5em, colsep = 2pt},
				cells = {m,c},
				hlines,
				vlines,
			}
			$\bullet$ &  \\
			$\bullet$ & $\bullet$ \\
		\end{tblr}
		\begin{tblr}{
				rows = {1.5em, rowsep = 2pt},
				columns = {1.5em, colsep = 2pt},
				cells = {m,c},
				hlines,
				vlines,
			}
			& $\bullet$ \\
			$\bullet$ &  \\
		\end{tblr}
		\caption{All $[2]^1$ and $[2]^2$ Local Arrays of Tropical Linear Series} 
		\label{Fig:rk2-d12-local-arrays}
	\end{figure}
	
	Local arrays of tropical linear series satisfy a number of combinatorial properties. Three properties correspond to closure under tropical linear combinations, realization of slopes on a tangent vector, and tropical dependence of functions.
	
	\begin{enumerate}
		\item [\namedlabel{prop:P1}{(P1)}] The set of dotted points is closed under pointwise minimum, i.e., the meet operation of the poset.
		\item [\namedlabel{prop:P2}{(P2)}] For all $1 \leq i \leq d$ and for all $1 \leq j \leq r+1$ there exists an element $\mathbf{w} \in \overline{P_v}$ such that $w_{i} = j$. That is, there exists a dotted point which is equal to $j$ in the $i^{\text{th}}$ dimension.
		\item [\namedlabel{prop:P3}{(P3)}] Any set $S$ of $r+2$ elements of $\overline{P_v}$ must contain a subset $S' \subseteq S$ such that for each $1 \leq i \leq d$, we have $\min \{ x_i : \mathbf{x} \in S' \}$ occurs at least twice.
	\end{enumerate}
	
	By slight abuse of notation, we will sometimes say that a set $S'$ of elements satisfies \ref{prop:P3} for $i$. This means that $\min \{ x_{i} : \mathbf{x} \in S' \}$ occurs at least twice, but that this may not be true for $j \neq i$.
	
	\begin{lemma}
		\label{Lem:tls-local-arrays}
		The local array $\overline{P_v}$ of a tropical linear series $(D,\Sigma)$ at a point $v \in \Gamma$ satisfies \ref{prop:P1}, \ref{prop:P2}, and \ref{prop:P3}.
	\end{lemma}
	
	\begin{proof}
		Let $\mathbf{x} = (x_1, \ldots, x_d)$ and $\mathbf{y} = (y_1, \ldots, y_d)$ be two elements of $\overline{P_v}$ and let $f, g \in \Sigma$ be associated to $\mathbf{x}$ and $\mathbf{y}$ respectively. Let $h \in \Sigma$ be $h = \min \{ f - f(v), g - g(v) \}$. We have $$\mathrm{sl}_{\eta_i}(h) = \min \{ \mathrm{sl}_{\eta_i}[x_i], \mathrm{sl}_{\eta_i}[y_i] \} = \mathrm{sl}_{\eta_i}[\min \{ x_i,y_i \}].$$ So, the associated element of $h$ in $\overline{P_v}$ is $\mathbf{x} \wedge \mathbf{y}$. This proves \ref{prop:P1}.
		
		For a tangent vector $\eta_i$, the slope $\mathrm{sl}_{\eta_i}[j]$ must be realized by some $f \in \Sigma$. So, the element of $\overline{P_v}$ associated to $f$ will have $i^{\text{th}}$ coordinate equal to $j$. This proves \ref{prop:P2}.
		
		Let $S$ be a set of $r+2$ elements of $\overline{P_v}$. There exists a set $T$ of functions $f_1, \ldots, f_{r+2} \in \Sigma$ associated to $S$. By Definition \ref{Def:Tropical-Linear-Series}, these functions must be tropically dependent. Let $T'$ be the set of functions attaining the minimum at $v$ in this tropical dependence.
		
		For each tangent vector $\eta_i \in T_v(\Gamma)$, there exist two functions $f,g \in T'$ such that $\mathrm{sl}_{\eta_i}(f) = \mathrm{sl}_{\eta_i}(g) \leq \mathrm{sl}_{\eta_i}(h)$ for all $h \in T'$. The set $S'$ of elements of $\overline{P_v}$ associated to $T'$ is a subset of $S$ with the required property for \ref{prop:P3}.
	\end{proof}
	
	We provide an example and non-example to further explore \ref{prop:P3}.
	
	\begin{example}
		For an example, consider the leftmost dot array of Figure \ref{Fig:ex-and-nonex-loc-arrays}. The four elements of this dot array contain the three dots in the top-right of the array (these three dots are the subset $S'$), which satisfy the necessary conditions. 
		
		As a non-example, the middle dot array of Figure \ref{Fig:ex-and-nonex-loc-arrays} fails to satisfy \ref{prop:P3}. Suppose an appropriate $S'$ exists. For the set $S'$ to satisfy \ref{prop:P3} for $i=1$ (the rows) it must contain the points $(3,1)$ and $(3,2)$ so that two elements are equal along this dimension. However, this means that $S'$ cannot satisfy \ref{prop:P3} for $i=2$ (the columns), as $(3,1)$ uniquely obtains the minimum along this dimension.
	\end{example}
	
	\begin{figure}
		\centering
		\begin{tblr}{
				rows = {1.5em, rowsep = 2pt},
				columns = {1.5em, colsep = 2pt},
				cells = {m,c},
				hlines,
				vlines,
			}
			& $\bullet$ &  \\
			& $\bullet$ & $\bullet$ \\
			$\bullet$ &  &  \\
		\end{tblr}
		\begin{tblr}{
				rows = {1.5em, rowsep = 2pt},
				columns = {1.5em, colsep = 2pt},
				cells = {m,c},
				hlines,
				vlines,
			}
			$\bullet$ & $\bullet$ &  \\
			& $\bullet$ &  \\
			& $\bullet$ &  \\
		\end{tblr}
		\begin{tblr}{
				rows = {1.5em, rowsep = 2pt},
				columns = {1.5em, colsep = 2pt},
				cells = {m,c},
				hlines,
				vlines,
			}
			$\bullet$ &  &  \\
			$\bullet$ &  &  \\
			$\bullet$ & $\bullet$ & $\bullet$ \\
		\end{tblr}
		\caption{Example and Non-Examples of $[3]^2$ Local Arrays of Tropical Linear Series}
		\label{Fig:ex-and-nonex-loc-arrays}
	\end{figure}
	
	Combinations of these properties allow us to prove additional facts about local arrays of tropical linear series.
	
	\begin{lemma}
		\label{Lem:Similar-Dotted-Entries}
		Let $P$ be a dot array of $[r+1]^d$ satisfying \ref{prop:P2} and \ref{prop:P3}. If $\mathbf{x}, \mathbf{y}$ are distinct elements of $P$ with $x_i = y_i$ for some $i$, then $x_i = y_i < r+1$.
	\end{lemma}
	
	\begin{proof}
		By \ref{prop:P2}, there exist $\mathbf{w}_1, \ldots, \mathbf{w}_r$ that are elements of $P$, such that $w_{j,i}$ are all distinct from each other and distinct from $x_i = y_i$. We now have a set $S = \{ \mathbf{x}, \mathbf{y}, \mathbf{w}_1, \dots, \mathbf{w}_r \}$ that must satisfy \ref{prop:P3}. We argue by constructing the set $S' \subseteq S$ that satisfies \ref{prop:P3} for all coordinates.
		
		The $i^{th}$ coordinate of two of the elements of $S$ must be equal. We know that this is only true for $\mathbf{x}$ and $\mathbf{y}$, so $\mathbf{x}, \mathbf{y} \in S'$. Because $\mathbf{x}$ and $\mathbf{y}$ are distinct, for some $k \neq i$ the two elements of $S'$ satisfying \ref{prop:P3} for $k$ cannot be $\mathbf{x}$ and $\mathbf{y}$, so there is some $1 \leq j \leq r$ such that $\mathbf{w}_j \in S'$. Therefore, $x_i = y_i \leq w_{j,i}$ and because these must be distinct, $x_i = y_i < w_{j,i}$, and hence $x_i = y_i < r+1$.
	\end{proof}
	
	\begin{lemma}
		\label{Lem:local-arrays-are-redundant-closures-low-rank}
		If $r = 0$ or $r = 1$ and $\overline{P} \subseteq [r+1]^d$ satisfies \ref{prop:P1}, \ref{prop:P2}, and \ref{prop:P3}, then $\overline{P}$ is the redundant closure of a permutation array of rank $r+1$ and dimension $d$.
		
		That is, the local array of a rank $r=0$ and $r=1$ tropical linear series at a point $v \in \Gamma$ is the redundant closure of a permutation array of rank $r+1$ and dimension $\mathrm{val}(v)$.
	\end{lemma}
	
	\begin{proof}
		For $r = 0$, a $[1]^d$-array has only one possible entry, and it must be dotted because the local array satisfies \ref{prop:P2}.
		
		For $r = 1$, we prove that $\overline{P}$ satisfies the definition of a totally rankable array given in Lemma \ref{Lem:Totally-Rankable-Characterization}.
		
		Let $\mathbf{x}, \mathbf{y}$ be elements of $\overline{P}$ with $x_i > y_i$ and $x_j = y_j$. By Lemma \ref{Lem:Similar-Dotted-Entries}, $x_j = y_j = 1 < 2$. By \ref{prop:P2}, there exists $\mathbf{z} \in \overline{P}$ such that $z_j = 2$. Consider $S = \{ \mathbf{x}, \mathbf{y}, \mathbf{z} \}$, which is a set of $r+2$ elements. To satisfy \ref{prop:P3} for $i$, $z_i = 1 = y_i < x_i$. For any other index $k$, $x_k = y_k$ implies that $x_k = y_k = 1$ by Lemma \ref{Lem:Similar-Dotted-Entries} and $x_k \neq y_k$ implies that $\min \{ x_k, y_k \} = 1$. In either case, $z_k \geq \min \{ x_k, y_k \}$.
		
		By Lemma \ref{Lem:Totally-Rankable-Arrays-Have-A-Permutation-Array} $\overline{P}$ contains a unique permutation array. By \ref{prop:P1}, $\overline{P}$ is the redundant closure of that permutation array.
	\end{proof}
	
	In Section \ref{Sec:Realizability} we will use \ref{prop:P3} when working with sparse permutation arrays. We will need the following result.
	
	\begin{lemma}
		\label{Lem:P3-and-perm-implies-redundant-perm}
		If a dot array $P$ satisfies \ref{prop:P1}, \ref{prop:P3}, and contains a sparse permutation array $Q$, then $P = \overline{Q}$.
	\end{lemma}
	
	\begin{proof}
		Because $P$ contains $Q$ and satisfies \ref{prop:P1}, it must also contain $\overline{Q}$. Let $\mathbf{q} \in \overline{Q} \backslash Q$. The set $S = Q \cup \{ \mathbf{q} \}$ is a set of $r+2$ elements of $P$, and thus must satisfy \ref{prop:P3}. As $Q$ is a sparse permutation array, all of its elements have unique values in every coordinate. Therefore, there exists a subset $S' \subseteq S$ such that $q_i = \min \{ x_i : \mathbf{x} \in S' \}$ for all $i$, and thus $\mathbf{q} \in \overline{Q}$.
	\end{proof}
	
	We now impose recursive conditions on $(D,\Sigma)$. We will primarily consider when $\Sigma$ is locally weakly recursive, but the following lemma will also consider an instance when $\Sigma$ is strongly recursive. First, we will need to expand on our understanding of subarrays.
	
	A $[k+1]^{d}$-\textit{subarray} of $[r+1]^d$ is $A_1 \times A_2 \times \cdots \times A_{d}$ where for all $i$, $A_i \subset [r+1]$ and $|A_i| = k+1$. Likewise, a $[k+1]^{d}$-\textit{subarray} of $P$ is the set of points $P \cap (A_1 \times \cdots \times A_{d})$ for some $[k+1]^{d}$-subarray of $[r+1]^d$.
	
	It is more natural to consider a $[k+1]^d$-subarray as actually a subset of $[k+1]^d$. To do so, we will use the following construction. Let $A_i = \{ a_{i,1} < \cdots < a_{i,k+1} \}$ and define the following functions: $$c_i : A_i \to [k+1]; a_{i,j} \mapsto j$$ and $$c : A_1 \times \cdots \times A_d \to [k+1]^d; (a_1, \ldots, a_d) \mapsto (c_1(a_1), \ldots , c_d(a_d)).$$ For a dot array $P \subseteq [r+1]^d$ and a $[k+1]^d$-subarray $P' = P \cap (A_1 \times \cdots \times A_d)$, we call $c(P')$ the \textit{shifted} $[k+1]^d$-subarray.
	
	\begin{enumerate}
		\item [\namedlabel{prop:P4}{(P4)}] For any set $S$ of $2$ elements of $\overline{P_v}$, there is a set $\overline{P_v}'$ and a $[2]^d$-subarray $(A_1 \times \cdots \times A_d) \cap \overline{P_v}$ such that $S \subseteq \overline{P_v}' \subseteq (A_1 \times \cdots \times A_d) \cap \overline{P_v}$ and the shifted subarray $c(\overline{P_v}')$ satisfies \ref{prop:P2} and \ref{prop:P3}.
		\item [\namedlabel{prop:P4-strong}{(P4-S)}] For any set $S$ of $1 \leq k+1 \leq r+1$ elements of $\overline{P_v}$, there is a set $\overline{P_v}'$ and a $[k+1]^d$-subarray $(A_1 \times \cdots \times A_d) \cap \overline{P_v}$ such that $S \subseteq \overline{P_v}' \subseteq (A_1 \times \cdots \times A_d) \cap \overline{P_v}$ and the shifted subarray $c(\overline{P_v}')$ satisfies \ref{prop:P2}, \ref{prop:P3}, and \ref{prop:P4-strong}.
	\end{enumerate}
	
	\begin{lemma}
		\label{Lem:recursive-local-arrays}
		The local array $\overline{P_v}$ of a rank $r$ tropical linear series $\Sigma$ at $v \in \Gamma$ satisfies \ref{prop:P4} if $\Sigma$ is locally weakly recursive and \ref{prop:P4-strong} if $\Sigma$ is strongly recursive.
	\end{lemma}
	
	\begin{proof}
		We prove the result when $\Sigma$ is strongly recursive. When $\Sigma$ is weakly recursive, the result is a special case of the strongly recursive condition. When $\Sigma$ is locally weakly recursive, we restrict to the subgraph $\mathrm{Star}_{\epsilon}(v)$ on which $\Sigma|_{\mathrm{Star}_{\epsilon}(v)}$ is weakly recursive.
		
		The result is vacuous when $r = 0$ and trivial when $r = 1$. Assume for induction that the statement is true for strongly recursive tropical linear series of rank $k < r$.
		
		Choose $1 \leq k+1 \leq r+1$ distinct elements of $\overline{P_v}$, and call this set $R$. There exist functions $f_1, \ldots, f_{k+1} \in \Sigma$ associated to the elements of $R$. These functions must be contained in a strongly recursive tropical linear subseries $\Sigma'$ of rank $k$. Let $\overline{P_v}' \subset \overline{P_v}$ be the elements associated to the functions of $\Sigma'$. Then, $R \subseteq \overline{P_v}' \subset \overline{P_v}$. We consider the shifted subarray $c(\overline{P_v}')$, which is the local array of $\Sigma'$ at $v$. By the inductive hypothesis, $c(\overline{P_v}') \subseteq [k+1]^d$ satisfies \ref{prop:P2}, \ref{prop:P3}, and \ref{prop:P4-strong}.
	\end{proof}
	
	We provide an example and non-example to further explore \ref{prop:P4}.
	
	\begin{example}
		To understand the recursive condition for non-trivial cases, we first must identify all of the $[2]^2$ local arrays of tropical linear series. By Lemma \ref{Lem:local-arrays-are-redundant-closures-low-rank} there are only two, pictured in Figure \ref{Fig:rk2-d12-local-arrays}. As mentioned in Subsection \ref{Subsec:Permutation-Arrays}, the dimension $2$ permutation arrays correspond to the nonzero elements of permutation matrices.
		
		Now we consider the leftmost and rightmost dot arrays in Figure \ref{Fig:ex-and-nonex-loc-arrays}. By checking every pair, one can see that the dot array on the left satisfies \ref{prop:P4}, as any choice of two dotted points is contained in a $[2]^2$-subarray whose dotted points can be shifted to resemble those in one of the arrays in Figure \ref{Fig:rk2-d12-local-arrays}.
		
		The dot array on the right of Figure \ref{Fig:ex-and-nonex-loc-arrays} does not satisfy \ref{prop:P4}. No $[2]^2$-subarray containing the top two rows contains dots that may be shifted to resemble one of the arrays of Figure \ref{Fig:rk2-d12-local-arrays}, so the elements $(1,2)$ and $(1,3)$ cannot be contained in an appropriate subarray.
	\end{example}
	
	With our light recursive conditions we may extend Lemma \ref{Lem:local-arrays-are-redundant-closures-low-rank}.
	
	\begin{theorem}
		\label{Thm:local-arrays-are-redundant-closures}
		The local array $\overline{P_v}$ of a rank $r$ locally weakly recursive tropical linear series $\Sigma$ at a point $v \in \Gamma$ is the redundant closure of a permutation array of rank $r+1$ and dimension $\mathrm{val}(v)$.
	\end{theorem}
	
	\begin{proof}
		For $r = 0$ and $r = 1$, the result was already proven in Lemma \ref{Lem:local-arrays-are-redundant-closures-low-rank}.
		
		For $r > 1$ we show that the argument reduces to the $r = 1$ case. Let $\mathbf{x}, \mathbf{y} \in \overline{P_v}$ with $x_i > y_i$ and $x_j = y_j$. By \ref{prop:P4}, $\mathbf{x}$ and $\mathbf{y}$ are contained in $\overline{P_v}' \subset \overline{P_v}$, which is a $[2]^d$-subarray such that $c(\overline{P_v}')$ satisfies \ref{prop:P2} and \ref{prop:P3}. We also know that $c(\overline{P_v}')$ satisfies \ref{prop:P1} because $\overline{P_v}$ satisfied $\ref{prop:P1}$.
		
		By Lemma \ref{Lem:local-arrays-are-redundant-closures-low-rank}, there exists $\mathbf{z}^c \in c(\overline{P_v}')$ such that $z_j^c = 2 > c_j(x_j) = c_j(y_j)$, $z_i^c = 1 = c_i(y_i)$, and $z_k^c \geq \min \{ c_k(x_k), c_k(y_k) \}$ for all $k \in \{ 1, \ldots , d \}$. Therefore, the element $\mathbf{z} = c^{-1}(\mathbf{z}^c) \in \overline{P_v}'$ satisfies $z_j > x_j = y_j$, $z_i = y_i$, and $z_k \geq \min \{ x_k, y_k \}$ for all $k \in \{ 1, \ldots , d \}$.
		
		By Lemma \ref{Lem:Totally-Rankable-Arrays-Have-A-Permutation-Array} $\overline{P_v}$ contains a unique permutation array. By \ref{prop:P1}, $\overline{P_v}$ is the redundant closure of that permutation array.
	\end{proof}
	
	In Section \ref{Sec:Realizability}, we investigate the converse of this theorem, i.e., when we can realize the redundant closure of a permutation array as the local array of a tropical linear series.
	
	At all but finitely many points we may further describe the local array of a locally weakly recursive tropical linear series. We write this result in slightly greater generality for later use.
	
	\begin{corollary}
		\label{Cor:most-arrays-standard-dim2}
		Let $\overline{P_v} \subseteq [r+1]^d$ be the local array of a tropical submodule $\Sigma \subseteq R(D)$ for a valence-2 point $v \in \Gamma$. Let $v \not \in \mathrm{supp}(D)$; let $\eta_1, \eta_2$ be the two tangent vectors of $v$; and let $\mathrm{sl}_{\eta_1}[i] + \mathrm{sl}_{\eta_2}[r+1-i] = 0$ for all $i$. If $\overline{P_v}$ is the redundant closure of a permutation array, then it is the redundant closure of the dimension $2$ standard permutation array.
	\end{corollary}
	
	\begin{proof}
		As $v \not \in \mathrm{supp}(D)$ and $\Sigma \subseteq R(D)$, for all $f \in \Sigma$ it must be that 
		
		\begin{align}
			\mathrm{sl}_{\eta_1}(f) + \mathrm{sl}_{\eta_2}(f) \leq 0. \label{Inequal:slope-restriction}
		\end{align}
		
		Let $f \in \Sigma$ be associated to the point $(i,x)$ for some $0 \leq i, x \leq r+1$. As $\mathrm{sl}_{\eta_1}[i] + \mathrm{sl}_{\eta_2}[r+1-i] = 0$, to satisfy (\ref{Inequal:slope-restriction}) we must have $x \leq r+1-i$. The only permutation matrix with all elements satisfying this condition is the anti-diagonal matrix, and therefore the local array must be the redundant closure of the standard permutation array.
	\end{proof}
	
	We finally arrive at the proof of the desired result.
	
	\begin{proof}[Proof of Theorem \ref{Thm:lwrtls-are-structured}]
		For every $v \in \Gamma$ and $\eta \in T_v(\Gamma)$ we have the appropriate collection of integers $S^{\eta}$ by Lemma \ref{Lem:r+1-slopes}. 
		
		By Theorem \ref{Thm:local-arrays-are-redundant-closures}, $\overline{P_v}$, the local array of $\Sigma$ at $v$, is the redundant closure of a rank $r+1$ and dimension $\mathrm{val}(v)$ permutation array, so we may choose $P_v$ to be the unique permutation array contained in $\overline{P_v}$.
		
		This defines an $(r+1)$-slope structure $\mathfrak{S}$, and $\Sigma \subseteq R(\mathfrak{S})$ by definition. Further, $\Sigma$ realizes $\mathfrak{S}$ as the element $\mathbf{x} \in \overline{P_v}$ is realized by some $f \in \Sigma$ by definition. 
		
		It remains only to show that the set $V$ of points that are not valence-2 or have an array other than the standard permutation array is finite, as all functions in the locally weakly recursive tropical linear series are compatible with the slope structure by construction. Let $V$ be the finite set supplied by Lemma \ref{Lem:Finite-slope-vectors}. Any point of $\Gamma \backslash V$ is valence-2, not in the support of $D$, and its tangent vectors $\eta_1$ and $\eta_2$ have the property $\mathrm{sl}_{\eta_1}[i] + \mathrm{sl}_{\eta_2}[r+1-i] = 0$. Corollary \ref{Cor:most-arrays-standard-dim2} implies that its local array is the redundant closure of the standard rank $r+1$ and dimension $2$ permutation array.
	\end{proof}
	
	\section{Extension of Theorem \ref{Thm:lwrtls-are-structured} and Counterexamples to the Converse of Corollary \ref{Cor:SRTLS-are-structured}}
	\label{Sec:extension-and-counterexample}
	
	We now extend Theorem \ref{Thm:lwrtls-are-structured} to all tropical linear series of small rank or on graphs with only points of small valence. We also provide an infinite class of counterexamples showing that not all structured tropical linear series are strongly recursive.
	
	\subsection{Extending Theorem \ref{Thm:lwrtls-are-structured}} 
	\label{Subsec:result-extension}
	
	Section \ref{Sec:LWRTLS-are-structured} leaves us with the following natural question.
	
	\begin{question}
		\label{Quest:tls-all-structured}
		Are all tropical linear series structured?
	\end{question}
	
	At present, this question remains open. One possible approach would be to prove that every tropical linear series satisfies a weaker recursive condition and then apply Theorem \ref{Thm:lwrtls-are-structured}. This motivates the following question, which is also of independent interest.
	
	\begin{question}
		\label{Quest:tls-are-locally-weakly-recursive}
		Are all tropical linear series weakly recursive? What about locally weakly recursive?
	\end{question}
	
	Another approach would be to combinatorially prove that local arrays of tropical linear series are redundant closures of permutation arrays using properties similar to \ref{prop:P1}, \ref{prop:P2}, and \ref{prop:P3} from Section \ref{Sec:LWRTLS-are-structured}. This is true for local arrays of dimension $d=1$ by \ref{prop:P2}, and has already been proven true for tropical linear series of rank $r=0$ and $r=1$ in Lemma \ref{Lem:local-arrays-are-redundant-closures-low-rank}. We extend our results slightly using this approach. We require a generalization of \ref{prop:P2}. 
	
	Recall the argument proving that a local array of a tropical linear series satisfies \ref{prop:P2}. Because a tropical linear series of rank $r$ has Baker-Norine rank $r$, it must contain a function that bends at any $r$ points of $\Gamma$. If we choose $r$ distinct points that are sufficiently close to a point $v$ all on the same edge $e_i$ that contains a tangent vector $\eta$ of $v$, this implies that there are at least $r+1$ slopes realized along this tangent vector. These $r+1$ slopes provide the necessary realization to prove \ref{prop:P2}. 
	
	The generalization of \ref{prop:P2} below takes the same approach, but considers choosing points along two distinct edges $e_i$ and $e_j$ containing tangent vectors $\eta_i$ and $\eta_j$ of $v$. We choose $r_i$ points on $e_i$ and $r_j$ points on $e_j$ such that $r_i + r_j = r$. This will imply the existence of points $\mathbf{t}, \mathbf{w}_1, \ldots, \mathbf{w}_{r_i}, \mathbf{u}_1, \ldots, \mathbf{u}_{r_j}$ in the local array satisfying the following: 
	
	\begin{itemize}
		\item In the $i^{\text{th}}$ dimension of the local array, the points $\mathbf{w}_1, \ldots, \mathbf{w}_{r_i}, \mathbf{t}$ form a strictly increasing chain.
		\item In the $j^{\text{th}}$ dimension of the local array, the points $\mathbf{u}_1, \ldots, \mathbf{u}_{r_j}, \mathbf{t}$ form a strictly increasing chain.
		\item In the $j^{\text{th}}$ dimension of the local array, any point $\mathbf{u}_{k}$ is strictly less than any point $\mathbf{w}_{k'}$.
	\end{itemize}
	
	This statement is made rigorous below.
	
	\begin{enumerate}
		\item [\namedlabel{prop:P2-G}{(P2-G)}] For any choice of fixed $i, j \in \{ 1, \ldots, d \}$ and for all partitions $r = r_{i} + r_{j}$ where $r_{i}, r_{j} \geq 0$, the set of dotted points of $\overline{P_v}$ contains points $\mathbf{t}, \mathbf{w}_1, \ldots, \mathbf{w}_{r_i}, \mathbf{u}_1, \ldots, \mathbf{u}_{r_j}$ such that the following are true:
		\begin{enumerate}
			\item [\namedlabel{prop:P2a}{(P2a)}] $w_{1,i} < w_{2,i} < \cdots < w_{r_{i}, i} < t_{i}$,
			\item [\namedlabel{prop:P2b}{(P2b)}] $u_{1,j} < u_{2,j} < \cdots < u_{r_{j}, j} < t_{j}$, and
			\item [\namedlabel{prop:P2c}{(P2c)}] $u_{k,j} < w_{k',j}$ for all $\mathbf{u}_k, \mathbf{w}_{k'}$.
		\end{enumerate}
	\end{enumerate}
	
	\begin{lemma}
		\label{Lem:tls-local-arrays-generalization}
		A local array $\overline{P_v}$ of a rank $r$ tropical linear series $(D, \Sigma)$ at a valence-$d$ point $v \in \Gamma$ satisfies \ref{prop:P2-G}.
	\end{lemma}
	
	\begin{proof}
		If $d = 1$, $r_i = 0$, or $r_j = 0$, then this property reduces to \ref{prop:P2}. We now consider when $d \geq 2$ and $r_i, r_j \geq 1$.
		
		We restrict our attention to $\mathrm{Star}_{\epsilon}(v)$. We restrict the $\epsilon$ so that $\mathrm{Star}(v) \cap \mathrm{supp}(D) \subseteq \{ v \}$ and that no generator of $\Sigma$ bends on $\mathrm{Star}(v) \backslash \{ v \}$. We label the edges of $\mathrm{Star}(v)$ as $e_1, \ldots, e_d$. The assumptions mean that for any $\mathbf{x} \in \overline{P_v}$, there is a function in $f_{\mathbf{x}} \in \Sigma$ which does not bend on $\mathrm{Star}(v) \backslash \{ v \}$ such that $f$ realizes $\mathbf{x}$. The set $\{ f_{\mathbf{x}} | \mathbf{x} \in \overline{P_v} \}$ is a generating set for $\Sigma|_{\mathrm{Star}_{\epsilon}(v)}$, which is a tropical linear series of rank $r$.
		
		Choose distinct points $p_1, \ldots, p_{r_i}$ on the interior of edge $e_i$ such that 
		
		\begin{align}
			(\mathrm{sl}_{\eta_i}[r+1] - \mathrm{sl}_{\eta_i}[1]) \cdot \max_{1 \leq k \leq r_i} \{ d(v,p_k) \} = \delta < \epsilon \label{Inequal:edge-i-restriction}
		\end{align}
		
		where $d(v,p_k)$ is the distance from $v$ to $p_k$.
		
		Choose distinct points $q_1, \ldots , q_{r_j}$ on the interior of edge $e_j$ such that 
		
		\begin{align}
			\delta < \min_{1 \leq k \leq r_j} \{ d(v,q_k) \}. \label{Inequal:edge-j-restriction}
		\end{align}
		
		Let $E_1 = p_1 + \cdots + p_{r_i}$, $E_2 = q_1 + \cdots + q_{r_j}$, and $E = E_1 + E_2$. Because $\Sigma$ is a tropical linear series, there exists $f \in \Sigma$ such that $\mathrm{div}(f) + D \geq E$. Let $\partial_v(f) = \mathbf{t} \in \overline{P_v}$. We argue by describing the generators of $\Sigma$ necessary to construct $f$. We let $\overline{f}$ be the function in the generating set associated to $\mathbf{t}$.
		
		To satisfy $\mathrm{div}(f) + D - E_1 \geq 0$, the function $f$ will decrease slope at the $r_i$ distinct points $p_1, \ldots, p_{r_i}$. This implies the existence of generators $g_1, \ldots, g_{r_i} \in \Sigma$ such that $$\mathrm{sl}_{\eta_i}(g_1) < \cdots < \mathrm{sl}_{\eta_i}(g_{r_i}) < \mathrm{sl}_{\eta_i}(\overline{f}).$$ We assume that these generators are chosen with the appropriate tropical scaling to construct $f$. That is, $\mathrm{div}(\min \{ \overline{f}, g_1, \ldots, g_{r_i} \}) + D - E_1 \geq 0$ with no further tropical scaling.
		
		Let $\mathbf{w}_1, \ldots, \mathbf{w}_{r_i} \in \overline{P_v}$ be the points associated to $g_1, \ldots, g_{r_i}$, respectively. Then, $w_{1,i} < w_{2,i} < \cdots < w_{r_i, i} < t_i$, which is \ref{prop:P2a}.
		
		For the remaining properties, we will need two separate cases. Pick $g_{k'}$ such that $g_{k',j} = \min_{1 \leq k \leq r_i} \{ g_{k,j} \}$. If multiple functions satisfy this condition, choose $g_{k'}$ such that $g_{k'}(v)$ is also minimized.
		
		\textbf{Case A:} If $\mathrm{sl}_{\eta_j}(\overline{f}) \leq \mathrm{sl}_{\eta_j}(g_{k'})$, then $\overline{f}(v') < g_{k'}(v')$ for all $v'$ on edge $e_j$. So, $\overline{f}$ will be the only relevant function for the construction of $f$ on edge $e_j$ among the generators of $f$ we have already identified. To satisfy $\mathrm{div}(f) + D - E_2 \geq 0$, the function $f$ will decrease in slope at the $r_j$ distinct points $q_1, \ldots, q_{r_j}$. This implies the existence of generators $h_1, \ldots, h_{r_j} \in \Sigma$ such that $$\mathrm{sl}_{\eta_j}(h_1) < \cdots < \mathrm{sl}_{\eta_j}(h_{r_j}) < \mathrm{sl}_{\eta_j}(\overline{f}).$$
		
		Let $\mathbf{u}_1, \ldots, \mathbf{u}_{r_j} \in \overline{P_v}$ be the points associated to $h_1, \ldots, h_{r_j}$, respectively. Then, $u_{1,j} < u_{2,j} < \cdots < u_{r_j, j} < t_j$, which is \ref{prop:P2b}, and because we assumed $\mathrm{sl}_{\eta_j}(\overline{f}) \leq \mathrm{sl}_{\eta_j}(g_{k'})$, these points of $\overline{P_v}$ will also satisfy \ref{prop:P2c}.
		
		\textbf{Case B:} If $\mathrm{sl}_{\eta_j}(\overline{f}) > \mathrm{sl}_{\eta_j}(g_{k'})$, we consider the following. Using the construction for $f$, we have the following:
		
		\begin{align*}
			\overline{f}(p_{k'}) & \geq g_{k'}(p_{k'}) \\
			\overline{f}(v) + \mathrm{sl}_{\eta_i} (\overline{f}) \cdot d(v,p_{k'}) & \geq g_{k'}(v) + \mathrm{sl}_{\eta_i} (g_{k'}) \cdot d(v,p_{k'}) \\
			\overline{f}(v) + (\mathrm{sl}_{\eta_i} (\overline{f}) - \mathrm{sl}_{\eta_i} (g_{k'})) \cdot d(v,p_{k'}) & \geq g_{k'}(v) \\
			\overline{f}(v) + \delta & \geq g_{k'}(v).
		\end{align*}
		
		\noindent The last inequality comes from the choice made in (\ref{Inequal:edge-i-restriction}).
		
		Let $v'$ be the point on $e_j$ such that $d(v,v') = \delta$. Let $\xi$ be the tangent vector of $v'$ directed away from $v$. Given that $\overline{f}(v) + \delta \geq g_{k'}(v)$, on the interval $[v,v']$ the slope of $\overline{f}$ is greater than the slope of $g_{k'}$, and $d(v,v') = \delta$, $\overline{f}(v') \geq g_{k'}(v')$. In fact, because $g_{k}(v) \geq \overline{f}(v)$ for all $k$ and our choice of $g_{k'}$, we also have that $g_{k}(v') \geq g_{k'}(v')$ for all $k$.
		
		Because of the choice made in (\ref{Inequal:edge-j-restriction}) and $g_{k'}$ is the minimum function at $v'$ among the generators we have identified, it the only function we need to consider as we construct $f$ along $e_j$. To satisfy $\mathrm{div}(f) + D - E_2 \geq 0$, the function $f$ will decrease in slope at the $r_j$ distinct points $q_1, \ldots, q_{r_j}$. This implies the existence of generators $h_1, \ldots, h_{r_j} \in \Sigma$ such that $$\mathrm{sl}_{\eta_j}(h_1) < \cdots < \mathrm{sl}_{\eta_j}(h_{r_j}) < \mathrm{sl}_{\eta_j}(g_{k'}).$$
		
		Let $\mathbf{u}_1, \ldots, \mathbf{u}_{r_j} \in \overline{P_v}$ be the points associated to $h_1, \ldots, h_{r_j}$, respectively. Then, $u_{1,j} < u_{2,j} < \cdots < u_{r_j, j} < w_{k',j}$. Because $\mathrm{sl}_{\eta_j}(\overline{f}) > \mathrm{sl}_{\eta_j}(g_{k'})$, \ref{prop:P2b} is true, and the choice of $g_{k'}$ implies that \ref{prop:P2c} is true.
	\end{proof}
	
	We can use this result to show that all tropical linear series are structured for low rank.
	
	\begin{corollary}
		\label{Cor:tls-structured-for-small-rank}
		If $(D, \Sigma)$ is a tropical linear series of rank $r = 0,1,2$ then it is structured.
	\end{corollary}
	
	\begin{proof}
		Once we show that all local arrays are the redundant closure of permutation arrays we may apply the argument as before to construct an appropriate slope structure. For $r = 0$ and $r = 1$ we have already proven this fact in Lemma \ref{Lem:local-arrays-are-redundant-closures-low-rank}. It remains to show for $r = 2$.
		
		Let $\mathbf{x}, \mathbf{y} \in \overline{P_v}$ with $x_i > y_i$ and $x_j = y_j$. By Lemma \ref{Lem:Similar-Dotted-Entries}, $x_j = y_j \neq 3$.
		
		If $x_j = y_j = 2$, then let $r_i = 0$ and $r_j = 2$. By \ref{prop:P2-G}, there are points $\mathbf{t}, \mathbf{u}_1, \mathbf{u}_{2} \in \overline{P}$ satisfying \ref{prop:P2b}. Because the local array is a subset of $[3]^d$ and \ref{prop:P2b} is satisfied, we know that $u_{k,j} = k$ and $t_j = 3$. 
		
		The set $S = \{ \mathbf{u}_1, \mathbf{x}, \mathbf{y}, \mathbf{t} \}$ must satisfy \ref{prop:P3}. The set $S' \subseteq S$ satisfying \ref{prop:P3} cannot include $\mathbf{u}_1$, as this element uniquely attains the minimum in the $j^{\text{th}}$-coordinate. Therefore, the set $S' = \{ \mathbf{x}, \mathbf{y}, \mathbf{t} \}$ must satisfy \ref{prop:P3}. This means that $t_i = y_i$, $t_k \geq \min \{ x_k, y_k \}$ for all $k$, and we already found that $t_j > x_j = y_j$.
		
		If $x_j = y_j = 1$, then let $r_i = 1$ and $r_j = 1$. By \ref{prop:P2-G}, there are points $\mathbf{t}, \mathbf{u}_1, \mathbf{w}_{1} \in \overline{P}$ satisfying \ref{prop:P2a}, \ref{prop:P2b}, and \ref{prop:P2c}. 
		
		The set $S = \{ \mathbf{x}, \mathbf{y}, \mathbf{t}, \mathbf{w}_1 \}$ must satisfy \ref{prop:P3}. That is, there is a set $S' \subseteq S$ satisfying \ref{prop:P3}.
		
		We now consider two cases.
		
		If $t_i = y_i$, then by \ref{prop:P2a}, $w_{1,i} < t_i$ implies that $\mathbf{w} \not \in S'$. Therefore, $t_k \geq \min \{ x_k, y_k \}$ in order for $S'$ to satisfy \ref{prop:P3}. We also know that that $t_j > u_{1,j} \geq 1$. So, $\overline{P_v}$ must be totally rankable.
		
		Suppose $w_{1,i} = y_i$ and $w_{1,k} \geq \min \{ x_k, y_k \}$ for all $k$. Then $w_{1,j} > u_{1,j} \geq 1$ implies that $\overline{P_v}$ is totally rankable. If there exists $k$ such that $w_{1,k} < \min \{ x_k, y_k \}$, then because $S = \{ \mathbf{x}, \mathbf{y}, \mathbf{t}, \mathbf{w}_1 \}$ satisfies \ref{prop:P3}, $t_k = w_{1,k}$. Further, $1 = t_k = w_{1,k} < \min \{ x_k, y_k \} = 2$ by Lemma \ref{Lem:Similar-Dotted-Entries}. By \ref{prop:P2-G}, there must exist some $\mathbf{t}' \in \overline{P_v}$ such that $t_j', t_k' \geq 2$. We have shown that $\mathbf{t}'$ is not equal to any of $\mathbf{x}, \mathbf{y}, \mathbf{t}, \mathbf{w}_1$. The set $S_{\text{new}} = \{ \mathbf{x}, \mathbf{y}, \mathbf{t}, \mathbf{t}' \}$ must contain a subset $S_{\text{new}}'$ satisfying \ref{prop:P3}. Because $t_k < \min \{ x_k, y_k, t_k' \}$, $\mathbf{t} \not \in S_{\text{new}}'$. So, $S_{\text{new}}' = \{ \mathbf{x}, \mathbf{y}, \mathbf{t}' \}$ satisfies \ref{prop:P3} implies that $\overline{P_v}$ is totally rankable.
	\end{proof}
	
	We can also combine properties of local arrays to achieve new results. We now restrict to tropical linear series of rank $r \geq 3$.
	
	\begin{lemma}
		\label{Lem:tls-local-arrays-nearly-totally-rankable}
		Let $P$ be a dot array of $[r+1]^d$ with $r \geq 3$ satisfying \ref{prop:P2-G} and \ref{prop:P3}. If $\mathbf{x}, \mathbf{y} \in P$ with $x_i > y_i$ and $x_j = y_j$ then there exists an element $\mathbf{z} \in \overline{P_v}$ such that $z_i = y_i$ and $z_j > y_j$.
	\end{lemma}
	
	\begin{proof}
		Let $r_i = r - y_j$ and $r_j = y_j$. We know that $r_i \geq 0$ and $r_j \geq 1$ because $y_j < r+1$ by Lemma \ref{Lem:Similar-Dotted-Entries}. Because $P$ satisfies \ref{prop:P2-G}, there are points $\mathbf{t}, \mathbf{w}_1, \ldots, \mathbf{w}_{r_i}, \mathbf{u}_1, \ldots, \mathbf{u}_{r_j} \in P$ satisfying \ref{prop:P2a}, \ref{prop:P2b}, and \ref{prop:P2c}.
		
		All of the values $u_{k,j}$ are distinct from each other and less than $t_j$ and $w_{k',j}$ by \ref{prop:P2b} and \ref{prop:P2c}. Therefore, $t_j \geq r_j$ and $\min{w_{k',j}} \geq r_j$. If we modify our choices of $\mathbf{u}_{k}$ to satisfy $u_{k,j} = k$ for $1 \leq k \leq r_{j} - 1$, then the elements $\mathbf{t}, \mathbf{w}_1, \ldots, \mathbf{w}_{r_i}, \mathbf{u}_1, \ldots, \mathbf{u}_{r_j}$ still satisfy \ref{prop:P2a}, \ref{prop:P2b}, and \ref{prop:P2c}. 
		
		Consider the set $S = \{ \mathbf{x}, \mathbf{y}, \mathbf{t}, \mathbf{w}_1, \ldots, \mathbf{w}_{r_i}, \mathbf{u}_1, \ldots, \mathbf{u}_{r_j - 1} \}$. This is a set of $3 + (r - y_j - 1) + y_j = r+2$ elements. Because $\overline{P_v}$ satisfies \ref{prop:P3}, $S$ must contain a set $S'$ such that for each $1 \leq k \leq d$, we have $\min \{ a_k : \mathbf{a} \in S' \}$ occurs at least twice. We argue through construction of $S'$.
		
		Consider the set $R = \{ \mathbf{a} \in S : a_j < y_j \} = \{ \mathbf{u}_1, \ldots, \mathbf{u}_{r_j - 1} \}$. All of the elements of $R$ have distinct values in the $j^{\text{th}}$ coordinate and for all $a \in R$ and all $b \in S \backslash R$, $b_j > a_j$. Therefore, $R \cap S' = \emptyset$ in order to satisfy \ref{prop:P3} for $j$.
		
		Consider the set $T = \{ \mathbf{a} \in S : a_j > y_j \} = \{ \mathbf{t}, \mathbf{w}_1, \ldots, \mathbf{w}_{r_i} \}$. All of the elements of $T$ are distinct in the $i^{\text{th}}$ coordinate. Therefore, $\mathbf{x} \in S'$ or $\mathbf{y} \in S'$. The inclusion of either element in $S'$ implies that both must be elements of $S'$ in order to satisfy \ref{prop:P3} for $j$. 
		
		Therefore, $\min \{ a_j : j \in S' \} = y_j = x_j$, and $\{ \mathbf{x}, \mathbf{y} \} \subset S' \subseteq T \cup \{ \mathbf{x}, \mathbf{y} \}$. Because $\mathbf{y} \in S'$, $\min \{ a_i : \mathbf{a} \in S' \} \leq y_i < x_i$. As all of the elements of $T$ are distinct in the $i^{\text{th}}$ coordinate, the only way to satisfy \ref{prop:P3} for $i$ is for one of $t_i, w_{1,i}, \ldots, w_{r_i, i}$ to be equal to $y_i$. This will be the desired element $\mathbf{z}$.
	\end{proof}
	
	\begin{corollary}
		\label{Cor:tls-structured-at-small-dimension-points}
		Let $(D, \Sigma)$ be a rank $r$ tropical linear series on a metric graph $\Gamma$. If $v$ is a point of $\Gamma$ with $\mathrm{val}(v) \leq 2$, then the local array of $\Sigma$ at $v$ is the redundant closure of a rank $r+1$ permutation array. Moreover, at all but finitely many of these points, the local array is the redundant closure of the rank $r+1$ and dimension $2$ standard permutation array.
	\end{corollary}
	
	\begin{proof}
		For a point of valence $1$ it follows from \ref{prop:P2} that the local array is the redundant closure of a permutation array. 
		
		For a point of valence $2$, Lemma \ref{Lem:tls-local-arrays-nearly-totally-rankable} is sufficient to prove that the local array is totally rankable by the characterization given in Lemma \ref{Lem:Totally-Rankable-Characterization}, as $z_k \geq \min \{ x_k,y_k \}$ for all dimensions, as there are no other dimensions. Thus, the local array is the redundant closure of a permutation array by Lemma \ref{Lem:Totally-Rankable-Arrays-Have-A-Permutation-Array} and \ref{prop:P1}.
		
		At any point not in the support of $D$, the local array must be the redundant closure of the rank $r+1$ and dimension $2$ standard permutation array by Corollary \ref{Cor:most-arrays-standard-dim2}.
	\end{proof}
	
	\begin{corollary}
		\label{Cor:tls-structured-on-loop-and-interval}
		If $\Gamma$ is the interval metric graph or the loop metric graph and $(D, \Sigma)$ is a tropical linear series then it is structured.
	\end{corollary}
	
	\begin{proof}
		Corollary \ref{Cor:tls-structured-at-small-dimension-points} implies that we may construct a suitable slope structure on any metric graph with points only of valence $1$ or $2$. This includes interval and loop metric graphs.
	\end{proof}
	
	\subsection{Counterexamples to the Converse of Corollary \ref{Cor:SRTLS-are-structured}}
	\label{Subsec:not-all-structured-are-strongly-recursive}
	
	The converse of Corollary \ref{Cor:SRTLS-are-structured} is vacuously true for strongly recursive tropical linear series of rank $r = 0$ and trivially true for rank $r = 1$. We provide counterexamples to show that the converse is false for strongly recursive tropical linear series of rank $r \geq 3$.
	
	\begin{example}
		\label{Ex:vamos-matroid}
		Let $\Gamma$ be the interval metric graph of length $8$ with left endpoint $v$. In \cite[Theorem~6.3]{chang2025matroidal}, the authors show that any valuated matroid $\Delta \subseteq \overline{\mathbb{R}}^{d+1}$ defines a tropical linear series $\Sigma \subseteq R(dv)$. By Corollary \ref{Cor:tls-structured-on-loop-and-interval}, $\Sigma$ is a structured tropical linear series. In \cite[Example~6.6]{chang2025matroidal}, it is shown that when $\Delta$ is the V\'amos matroid $V$ this submodule does not satisfy condition (3) of \ref{Def:Strongly-Recursive-Linear-Series}, and is therefore not a strongly recursive tropical linear series.
	\end{example}
	
	This example originally appeared in \cite{chang2025matroidal} showing that not all matroidal linear series are strongly recursive, and thus not all matroidal linear series are tropicalizations of algebraic linear series. This is reflected with the transverse lines between these notions in Figure \ref{Fig:linear-series-poset}.
	
	\begin{remark}
		While discussing this example with Dave Jensen and Jidong Wang, we found it is possible to show that $\Sigma$ is also a weakly recursive tropical linear series. Hence, Figure \ref{Fig:linear-series-poset} includes two transverse lines between strongly recursive and weakly recursive tropical linear series. Further, this provides another proof that the example is a structured tropical linear series, as all weakly recursive tropical linear series are structured. 
		
		To see that $\Sigma$ is weakly recursive, let $f$ and $g$ be two functions of $\Sigma$. Let $p$ and $q$ be two elements in $\Delta$ that map to $f$ and $g$ respectively under the map $\Phi$ given in \cite[Theorem~6.3]{chang2025matroidal}. By \cite[Theorem~3.1.9]{wangPhDthesis}, these two points are contained in a tropical line that is also contained in $\Delta$, and by \cite[Theorem~6.3]{chang2025matroidal} the image of this tropical line under $\Phi$ is a tropical linear series of rank $1$.
	\end{remark}
	
	We have the following technical lemma, restated for the purposes of this paper.
	
	\begin{lemma}\cite[Proposition~7.5.5]{bjorner1999oriented}
		\label{Lem:LIP-extension}
		Let $M$ be a rank $r$ matroid on ground set $E$ that fails the Levi Intersection Property. A single element free extension $M'$ of $M$ is a rank $r+1$ matroid on ground set $E' = E \cup \{ e' \}$ that also fails the Levi Intersection Property.
	\end{lemma}
	
	By \cite[Theorem~6.3]{chang2025matroidal}, we may construct a tropical linear series of rank $r \geq 4$ on the interval, which by Corollary \ref{Cor:tls-structured-on-loop-and-interval} is structured. Using the construction in \cite[Example~6.6]{chang2025matroidal}, the constructed matroid will contain a set of $r$ flats that fail the Levi Intersection Property and thus also the Discrete Intersection Property. The set of $r$ flats corresponds to a set of $r$ functions of the tropical linear series which cannot be contained in a tropical linear subseries of rank $r-1$ without contradicting \cite[Corollary~4.8]{chang2025matroidal}. This provides the desired counterexample for $r \geq 4$.
	
	For more information regarding the construction of these counterexamples, the interested reader should consult \cite[Example~6.6]{chang2025matroidal} and the proof of \cite[Theorem~6.4]{wang2024lorentzian}.
	
	The following question remains open.
	
	\begin{question}
		\label{Quest:CLLS-are-SRTLS-r=2}
		Are all rank $2$ tropical linear series strongly recursive?
	\end{question}
	
	\section{Realizability of Permutation Arrays as Local Arrays of Tropical Linear Series}
	\label{Sec:Realizability}
	
	In this section we determine when the converse of Theorem \ref{Thm:local-arrays-are-redundant-closures} holds for strongly recursive tropical linear series, when it holds for arbitrary tropical linear series, and we provide an infinite class of counterexamples to realization by strongly recursive tropical linear series. That is, we consider when a permutation array can be realized as the local array of a tropical linear series.
	
	Throughout this section let $P$ be an $[r+1]^d$ permutation array and let $\Gamma$ be the star metric graph with central vertex $v$ with $d$ edges of length $\ell$ ordered $e_1, \ldots, e_d$. For $\mathbf{x} \in \overline{P}$, let $f_{\mathbf{x}}$ be the function with slope $x_i - 1$ along edge $e_i$. We investigate $\Sigma_P = \langle f_{\mathbf{x}} | \mathbf{x} \in \overline{P} \rangle \subseteq R(s \cdot v)$ for sufficiently large $s$.
	
	\subsection{Realizability by Strongly Recursive Tropical Linear Series}
	\label{Subsec:real-by-SRTLS}
	
	When $r = 0$ or $d = 1$, realizability as local arrays of strongly recursive tropical linear series of rank $r$ follows trivially from definitions.
	
	The complete linear series on a tree is the tropicalization of an algebraic linear series, and thus by \cite[Section~6]{farkas2025kodairadimensionsoverlinemathcalm22overlinemathcalm23} is a strongly recursive tropical linear series. So, the local array of $R(r \cdot v)$ at the central vertex of $\Gamma$ is the standard permutation array of rank $r+1$ and dimension $d$. Next, we have our first nontrivial realizability result.
	
	\begin{proposition}
		\label{Prop:realizability-of-sparse-perm-arrays}
		If $P$ is an $[r+1]^d$ sparse permutation matrix with $d \geq 2$ and $r \geq 1$, then there exists a rank $r$ strongly recursive tropical linear series with $\overline{P}$ as a local array of some point $p \in \Gamma$.
	\end{proposition}
	
	\begin{proof}
		Consider the rank $rd$ tropical linear series $R(rd \cdot v)$. Because $R(rd \cdot v)$ is a strongly recursive tropical linear series, the $r+1$ functions of the set $\{ f_{\mathbf{x}} | \mathbf{x} \in P \}$ must be contained in a strongly recursive tropical linear subseries $\Sigma \subseteq R(rd \cdot v)$.
		
		Consider the local array of $\Sigma$ at $v$. This local array contains the sparse permutation matrix $P$, and being a local array of a tropical linear series must satisfy properties \ref{prop:P1} and \ref{prop:P3}. Therefore, the local array of $\Sigma$ at $v$ must be exactly $\overline{P}$ by Lemma \ref{Lem:P3-and-perm-implies-redundant-perm}.
	\end{proof}
	
	This result implies realizability for $d=2$, as all arrays of this dimension are sparse.	Next, we prove realizability for permutation arrays of rank $r+1=2$ to be local arrays of strongly recursive tropical linear series. We do this by directly proving that $\Sigma_P$ is a tropical linear series of rank $r=1$. First, we need information and lemmas regarding projections of totally rankable arrays.
	
	Define the projection map $$\pi_i:[r+1]^d \to [r+1]^{d-1}, (x_1, \ldots, x_{i-1}, x_i, x_{i+1}, \ldots, x_d) \mapsto (x_1, \ldots, x_{i-1}, x_{i+1}, \ldots, x_d).$$
	
	\begin{lemma}
		\label{Lem:Proj-of-Perm}
		Let $P$ be an $[r+1]^d$ permutation array with $d \geq 2$ and let $S \subseteq P$. Then $|\pi_{i}(S)| = |S|$.
	\end{lemma}
	
	\begin{proof}
		Suppose $|\pi_{i}(S)| < |S|$. Then there exist distinct elements $\mathbf{x}, \mathbf{y} \in P$ such that $\pi_i(\mathbf{x}) = \pi_i(\mathbf{y})$. Because $d \geq 2$, $\mathbf{x}$ and $\mathbf{y}$ are equal in every coordinate except $i$, and $\mathbf{x} \wedge \mathbf{y} = \mathbf{x}$ or $\mathbf{x} \wedge \mathbf{y} = \mathbf{y}$. Therefore, one of $\mathbf{x}$ or $\mathbf{y}$ is redundant. This contradicts that $P$ is a permutation array.
	\end{proof}
	
	\begin{lemma}\cite[Lemma~3.4]{eriksson2000combinatorial}
		\label{Lem:proj-of-totally-rankable}
		If $P$ is totally rankable, then $\pi_i(P)$ is totally rankable.
	\end{lemma}
	
	To show that the generators of $\Sigma_P$ satisfy the tropical dependence condition, it suffices to show that the elements of $P$ satisfy \ref{prop:P3}, as all of the generators have constant slope on the edges of $\Gamma$. Then, $\Sigma_P$ satisfies the tropical dependence condition by Lemma \ref{Lem:generators-determine-dependence}.
	
	\begin{proposition}
		\label{Prop:realization-for-rk1-perm arrays}
		If $P$ is a $[2]^d$ permutation matrix then there exists a rank $1$ strongly recursive tropical linear series with $\overline{P}$ as a local array of some point $v \in \Gamma$.
	\end{proposition}
	
	\begin{proof}
		For $d=1$ and $d=2$ we already know that the result holds. Now assume $d \geq 3$.
		
		By definition, $\Sigma_P \subseteq R(d \cdot v)$ is a finitely generated tropical submodule. 
		
		As $P$ is rank $2$, for any $i \in \{ 1, \ldots, d \}$ there exists $\mathbf{x}, \mathbf{y} \in P$ such that $x_i = 1 < 2 = y_i$. Let $w$ be a point on the edge $e_i$. Then there exists $a \in [0, \ell]$ such that the function $\min \{ f_{\mathbf{y}}, f_{\mathbf{x}} + a \}$ satisfies \ref{prop:BNRP} for $E = w$.
		
		As mentioned above, it suffices to show that the elements of $P$ satisfy \ref{prop:P3}. If $P$ is sparse, then $|P| = 2$, so \ref{prop:P3} holds vacuously. If $P$ is not sparse, then $|P| \geq 3$. Consider $S \subseteq P$ such that $|S| = 3$.
		
		Assume for induction that the statement is true for $d-1 \geq 2$. Elements of $R(P) \backslash P$ can be generated by elements in $P$, so our inductive hypothesis implies that totally rankable arrays of dimension $d-1 \geq 2$ also satisfy \ref{prop:P3} by Lemma \ref{Lem:generators-determine-dependence}.
		
		Let $S_i = \pi_i^{-1}(\pi_i(S))$. We know that $\pi_1(S)$ and $\pi_2(S)$ are both subsets of totally rankable arrays of dimension $d-1$, so by the inductive hypothesis and Lemma \ref{Lem:proj-of-totally-rankable}, $S_1$ satisfies \ref{prop:P3} for $i \neq 1$ and $S_2$ satisfies \ref{prop:P3} for $i \neq 2$. 
		
		Let $S_i' \subseteq S_i$ be subset that satisfies \ref{prop:P3} for $j \neq i$. Because $d-1 \geq 2$, a set of two distinct elements of a projection cannot alone satisfy \ref{prop:P3}. Therefore, $|S_i'| = 3$ for $i=1,2$ and thus $S_i = S$. This means that $S$ satisfies \ref{prop:P3} for all $i$, and thus satisfies \ref{prop:P3}.
		
		The remaining conditions for strongly recursive tropical linear series are trivial when the series is rank $1$.
	\end{proof}
	
	\subsection{Realizability as Local Arrays of Tropical Linear Series}
	\label{Subsec:real-by-tls}
	
	We now work to extend realizability to higher rank, but only as local arrays of structured tropical linear series. That is, we will approach the following questions:
	
	\begin{question}
		Are all permutation arrays realizable as local arrays of structured tropical linear series? What about tropical linear series?
	\end{question}
	
	We leave open these questions, which are in fact equivalent. Say the local array of a tropical linear series $\Sigma$ at a point $v$ is $\overline{P}$, which is the redundant closure of a permutation array $P$. Consider the tropical linear series $\Sigma|_{\mathrm{Star}_{\epsilon}(v)} \subseteq R(s \cdot v)$ where $\epsilon$ is sufficiently small so that the generators of $\Sigma$ do not bend away from $v$ and $s$ is sufficiently large. On a star metric graph, we may construct a suitable slope structure for $\Sigma|_{\mathrm{Star}_{\epsilon}(v)}$ at all points away from the central vertex by Corollary \ref{Cor:tls-structured-at-small-dimension-points}. So, if the local array of the central point is the redundant closure of a permutation array, then the tropical linear series is structured. 
	
	We will answer the question for permutation arrays of rank $r+1=3$ and leave the general case open. It suffices to show that $\Sigma_P$ satisfies \ref{prop:BNRP} for all effective divisors of degree-$2$ and that the elements of $P$ satisfy \ref{prop:P3}. Because of the length of the proofs, we will split up this argument.
	
	\begin{lemma}
		\label{Lem:bn-realization-rk2}
		The tropical submodule $\Sigma_P$ satisfies \ref{prop:BNRP} for all effective divisors $E = v_1 + v_2$ on $\Gamma$.
	\end{lemma}
	
	\begin{proof}
		We already know that permutation arrays of dimension $1$ and $2$ can be realized as local arrays of tropical linear series, and thus $\Sigma_P$ satisfies \ref{prop:BNRP} for all effective divisors $E = v_1 + v_2$ on $\Gamma$ in these cases. Assume for induction that if $P$ is a rank $3$ and dimension $d-1 \geq 2$ permutation array, then $\Sigma_P$ satisfies \ref{prop:BNRP} for all degree-$2$ divisors.
		
		Let $P$ be a permutation array of dimension $d \geq 3$. Because $E$ is an effective divisor of degree $2$ and $d \geq 3$, there is always at least one edge $e_i$ of $\Gamma$ in which $\mathrm{supp}(E) \cap e_i \subseteq \{ v \}$, where $v$ is the central vertex. 
		
		By Lemma \ref{Lem:proj-of-totally-rankable}, $\pi_i(P)$ is a totally rankable array of rank $3$ and dimension $d-1$. By Lemma \ref{Lem:Totally-Rankable-Arrays-Have-A-Permutation-Array}, there is a permutation array $P' \subseteq \pi_i(P)$ which is rank $3$ and dimension $d-1$. By the inductive hypothesis, there exists a tropical combination $$f = \min_{\mathbf{x} \in P'} \{ f_{\mathbf{x}} + a_{\mathbf{x}} \}$$ which satisfies \ref{prop:BNRP} for $E$. By Lemma \ref{Lem:Proj-of-Perm}, each element $\mathbf{x} \in P'$ has a unique lift $\Tilde{\mathbf{x}} \in P$. So, the tropical combination $$\Tilde{f} = \min_{\mathbf{x} \in P'} \{ f_{\Tilde{\mathbf{x}}} + a_{\mathbf{x}} \}$$ will also satisfy \ref{prop:BNRP} for $E$, as on the $d-1$ edges containing the support of $E$, $f$ and $\Tilde{f}$ agree.
	\end{proof}
	
	The proof for tropical dependence for permutation arrays of rank $r+1=3$ readily uses the fact that the tropically dependent sets must be sufficiently small. This makes it easier to show that two projections must agree.
	
	\begin{lemma}
		\label{Lem:tropdep-realization-rk2}
		If $P$ is a $[3]^d$ permutation array then the elements of $P$ satisfy Property \ref{prop:P3}.
	\end{lemma}
	
	\begin{proof}
		Consider $S \subseteq P$ such that $|S| = 4$. Assume that the statement is true for $d-1 \geq 2$. This also means that totally rankable arrays of dimension $d-1$ satisfy Property \ref{prop:P3} by Lemma \ref{Lem:generators-determine-dependence}.
		
		Let $S_i = \pi_i^{-1}(\pi_i(S))$. We know that each $\pi_i(S)$ is a subset of a totally rankable array of dimension $d-1$, so by the inductive hypothesis and Lemma \ref{Lem:proj-of-totally-rankable}, $S_i$ satisfies \ref{prop:P3} for $j \neq i$. Let $S_i' \subseteq S_i$ be a subset that satisfies \ref{prop:P3} for $j \neq i$. Choose $S_i'$ such that $|S_i'|$ is maximal. 
		
		Because $d-1 \geq 2$, a set of two distinct elements of a projection cannot satisfy \ref{prop:P3}. Therefore, $|S_i'| \geq 3$.
		
		If there is $i_0$ such that $|S_{i_0}'| = 4$, then $S_{i_0}' = S$. If $S_{i_0}' = S$ passes \ref{prop:P3} directly, we are done. Otherwise, there exists $\mathbf{x} \in S$ such that $x_{i_{0}} < y_{i_0}$ for all $\mathbf{y} \in S \backslash \{ \mathbf{x} \}$. Therefore, $\mathbf{x} \not \in S_j'$ for $j \neq i_0$, since it uniquely attains the minimum on this coordinate among all elements of $S$. For $j, k, i_0$ all distinct we must have that $S_j', S_k' \subseteq S \backslash \{ \mathbf{x} \}$. Since $|S \backslash \{ \mathbf{x} \}| = 3$ and all $|S_i'| \geq 3$, we must have that $S_j' = S_k'$. Therefore, $S_j' = S_k'$ satisfies \ref{prop:P3} for all indices.
		
		If $|S_i'| = 3$ for all $i$. Because $|S_1'|$ is maximal, then for $\{ \mathbf{x} \} = S \backslash S_1'$ either $x_2 < y_2$ for all $\mathbf{y} \in S_1'$ or $x_3 < y_3$ for all $\mathbf{y} \in S_1'$. The former case implies that $\mathbf{x} \not \in S_3'$. Because $S_1', S_3' \subseteq S \backslash \{ \mathbf{x} \}$ and all $|S_i'| \geq 3$, we must have that $S_1' = S_3'$. Therefore, $S_1' = S_3'$ satisfies \ref{prop:P3} for all indices. A similar argument follows if $x_3 < y_3$ for all $\mathbf{y} \in S_1'$.
	\end{proof}
	
	Through a computer check, we have found that the tropical dependence property holds at least for permutation arrays of rank $r+1=4$ with $d \leq 4$ and $r+1 = 5$ with $d=3$.
	
	On a related note, Amini and Gierczak proved Lemma \ref{Lem:Admissible-implies-structured}, which is a realization result for slope structures by admissible submodules, under the assumption that a crude linear series of rank $r$ admits an admissible submodule of the same rank. However, they are unable to prove that this assumption always holds. The numerous cases involved in Lemma \ref{Lem:bn-realization-rk2} illustrate the difficulty of verifying this assumption for higher values of $r$ and more complex slope structures.
	
	\subsection{Counterexamples to Realizability as Local Arrays of Strongly Recursive Tropical Linear Series}
	\label{Subsec:counterexample-to-SRTLS-realizability}
	
	For higher rank, it is nontrivial to show that $\Sigma_P$ will satisfy the remaining conditions for a strongly recursive tropical linear series. Our counterexamples exploit this fact by providing a class of permutation arrays which do not satisfy \ref{prop:P4-strong}. The initial permutation array was originally given in \cite{billey2008intersections} as a counterexample to Eriksson and Linusson's conjecture regarding permutation arrays being realized by flag arrangements.
	
	\begin{example}\cite[Counterexample~3]{billey2008intersections}
		\label{Ex:counterexample-to-realizability}
		Consider the rank $r+1 = 4$ and dimension $d = 4$ permutation array $P$ given by the points below.
		
		\begin{center}
			$(1,3,1,4), (3,1,1,3), (1,1,2,3), (2,2,1,2), (2,1,2,2), (1,2,2,2),$
			$(1,1,3,2), (4,1,4,1), (3,4,1,1), (2,2,2,1), (1,3,3,1)$
		\end{center}
		
		To show that these points form a permutation array, one may use a computer check to verify that the set satisfies the condition given in Lemma \ref{Lem:Totally-Rankable-Characterization} and that none of the elements of this set are redundant points.
		
		Consider the three elements $A = \{ (3, 4, 1, 1), (1, 1, 3, 2), (3, 1, 1, 3) \}$. Again using a computer check, it is possible to show that there is no subset $\overline{P}'$ of the elements of $\overline{P}$ containing $A$ such that $\overline{P}'$ is a totally rankable $[3]^4$-subarray. Therefore, $P$ does not satisfy \ref{prop:P4-strong}, so it cannot be the local array of a strongly recursive tropical linear series.
	\end{example}
	
	We extend this counterexample to higher dimensions.
	
	\begin{lemma}
		There exist permutation arrays that do not satisfy \ref{prop:P4-strong} for $r+1 = 4$ and all $d \geq 4$.
		\label{Lem:counterexample-to-realizability-dimension-induction}
	\end{lemma}
	
	\begin{proof}
		We prove this by induction on $d$, with Example \ref{Ex:counterexample-to-realizability} as the base case. That is, assume that there exists $P$ that is a rank $4$, dimension $d-1 \geq 4$ permutation array that fails \ref{prop:P4-strong} via $A \subseteq P$ and $P$ is an antichain. One may verify that the permutation array in Example \ref{Ex:counterexample-to-realizability} is an antichain.
		
		Define $Q = \{ (x_1, \ldots, x_{d-1}, x_{d-1}) : (x_1, \ldots, x_{d-1}) \in P \}$. That is, the $d^{th}$ coordinate is equal to the $(d-1)^{th}$ coordinate. One may verify that $Q$ is a totally rankable array via the inductive hypothesis and Lemma \ref{Lem:Totally-Rankable-Characterization}. Further, $Q$ has no redundant points because $Q$, like $P$, is an antichain.
		
		Let $B = \{ (x_1, \ldots, x_{d}) : (x_1, \ldots, x_{d-1}) \in A \} \subseteq Q$. Suppose there exists $\overline{Q}'$ such that $B \subseteq \overline{Q}' \subseteq \overline{Q}$ and $\overline{Q}'$ is a totally rankable $[3]^d$-subarray. By Lemma \ref{Lem:proj-of-totally-rankable}, $A \subseteq \pi_{[d] \backslash \{d \}}(\overline{Q}') \subseteq P$ and $\pi_{[d] \backslash \{d \}}(\overline{Q}')$ must be a totally rankable $[3]^{d-1}$-subarray. However, this contradicts the inductive hypothesis. So, $Q$ fails \ref{prop:P4-strong} via $B$.
	\end{proof}
	
	We now show that we may extend our counterexample to all larger cases of rank and dimension.
	
	\begin{proposition}
		There exists a permutation array $P$ of rank $r+1$ and dimension $d$ that does not satisfy \ref{prop:P4-strong} for all $r+1 \geq 4$ and $d \geq 4$.
		\label{Prop:counterexample-to-realizability-rank-induction}
	\end{proposition}
	
	\begin{proof}
		We prove this by induction on $r+1$, with Lemma \ref{Lem:counterexample-to-realizability-dimension-induction} providing the base case for $r+1 = 4$. That is, assume that there exists $P$ that is a rank $r+1 \geq 4$, dimension $d$ permutation array that fails \ref{prop:P4-strong} via $A \subseteq P$.
		
		Define $Q = P \cup \{ (r+2, \ldots, r+2) \}$. One may verify that $Q$ is an $[r+2]^d$ permutation array via the inductive hypothesis and Lemma \ref{Lem:Totally-Rankable-Characterization}. The set $A$ as defined in Lemma \ref{Lem:counterexample-to-realizability-dimension-induction} has three distinct values in the last coordinate, so if $\overline{Q}'$ is a $[3]^d$-subarray of $\overline{Q}$ and $A \subseteq Q$, then $(r+2, \ldots, r+2) \not \in \overline{Q}'$, which implies that $\overline{Q}'$ is a $[3]^d$-subarray of $P$, which by the inductive hypothesis cannot be totally rankable. Thus, $Q$ fails \ref{prop:P4-strong} via $A \subseteq Q$ .
	\end{proof}
	
	These counterexamples leave open the following question:
	
	\begin{question}
		Are permutation arrays of rank $r+1=3$ realizable as local arrays of strongly recursive tropical linear series? What about of dimension $d=3$?
	\end{question}
	
	A computer check shows that all permutation arrays with rank $r+1 = 3$ and dimension $d \leq 5$ as well as those with dimension $d = 3$ and rank $r+1 \leq 5$ satisfy \ref{prop:P4-strong}. However, this does not imply that all necessary collections of functions in $\Sigma$ satisfy the recursion conditions, only the generators.
	
	\bibliography{math.bib}
	
\end{document}